\def\datum{ February 19, 2002}
\documentclass[11pt]{amsart}

\usepackage{amssymb}

\newcommand\al{\alpha}
\newcommand\bt{\beta}
\newcommand\G{\Gamma}
\newcommand\g{\gamma}
\newcommand\Dt{\Delta}
\newcommand\dt{\delta}
\newcommand\e{\varepsilon}
\renewcommand\th{\vartheta}
\renewcommand\k{\kappa}
\newcommand\Ld{\Lambda}
\newcommand\ld{\lambda}
\newcommand\m{\mu}
\newcommand\n{\nu}
\newcommand\x{\xi}
\newcommand\p{\pi}
\renewcommand\r{\rho}
\newcommand\s{\sigma}

\newcommand\ch{\chi}
\newcommand\ps{\psi}

\newcommand\z{\zeta}

\newcommand\irr{\varpi}
\newcommand\fka{{\mathfrak a}}
\newcommand\fkb{{\mathfrak b}}
\newcommand\fkc{{\mathfrak c}}
\newcommand\fkd{{\mathfrak d}}

\newcommand\fkp{{\mathfrak p}}
\newcommand\qq{{\mathfrak q}}

\newcommand\bhv{{\mathfrak H}}

\newcommand\GG{{\mathbf G}}

\newcommand\CC{\mathbb C}
\newcommand\QQ{\mathbb Q}
\newcommand\RR{\mathbb R}
\newcommand\ZZ{\mathbb Z}

\newcommand\Akr{{\mathcal A}}

\newcommand\Fkr{{\mathcal F}}
\newcommand\Hkr{{\mathcal H}}

\newcommand\Lkr{{\mathcal L}}
\newcommand\Mkr{{\mathcal M}}
\newcommand\Okr{{\mathcal O}}
\newcommand\Pkr{{\mathcal P}}
\newcommand\Rkr{{\mathcal R}}

\newcommand\Xkr{{\mathcal X}}
\newcommand\Ykr{{\mathcal Y}}

\newcommand\SL{{\mathrm{SL}}}
\newcommand\PSL{{\mathrm{PSL}}}
\newcommand\SO{{\mathrm{SO}}}

\newcommand\B{{\sf B}}

\newcommand\Eta{{\sf H}}
\newcommand\N{{\sf N}}

\newcommand\re{\operatorname{Re}}
\newcommand\im{\operatorname{Im}}
\newcommand\Tr{\operatorname{Tr}}
\newcommand\Gf{\operatorname{\Gamma}}
\newcommand\sign{\operatorname{sign}}
\newcommand\supp{\operatorname{Supp}}
\newcommand\oh{O}
\newcommand\vol{\operatorname{vol}}

\renewcommand\setminus{\smallsetminus}

\newcommand\one{{\bf 1}}

\makeatletter 
\newcommand\matc[4]{\left( {#1\@@atop #3}{#2\@@atop #4}\right)}
\newcommand\matr[4]{\left( {\hfill #1\@@atop\hfill #3}{\hfill
#2\@@atop\hfill #4}\right)}
\makeatother

\newcommand\widearray[1]{\renewcommand\arraystretch{1.4}
\begin{array}{#1}}

\newcommand\isdef{\mathrel{:\mskip2mu=}}

\newcommand\vz[1]{\mathchoice{\left\{ #1 \right\}}{\left\{ #1
\right\}}{\{ #1 \}}{\{ #1 \}}}

\newcommand\vzm[2]{\mathchoice{\left\{\, #1 : #2 \,\right\}}{\{\, #1
:\allowbreak #2 \,\}}{\{ #1 :\allowbreak #2 \}}{\{ #1 :\allowbreak #2
\}}}

\newcommand\Ndiscr{\operatorname{N^{\mathrm{discr}}}}

\theoremstyle{plain}

\newtheorem{thm}{Theorem}[section]
\newtheorem{lem}[thm]{Lemma}
\newtheorem{prop}[thm]{Proposition}
\newtheorem{cor}[thm]{Corollary}

\theoremstyle{definition}
\newtheorem{defn}[thm]{Definition}

\theoremstyle{remark}
\newtheorem{remark}[thm]{Remark}

\newcommand\rmk[1]{\medskip\par\noindent{\em #1. }\ignorespaces}

\begin{document}
\title[Density of automorphic forms on Hilbert modular groups]{Density
results for automorphic forms on Hilbert modular groups}

\author{R.W.\,Bruggeman}
\address{Mathematisch Instituut Universiteit Utrecht, Postbus 80010,
NL-3508 TA Utrecht, Nederland}
\email{bruggeman@math.uu.nl}

\author{R.J.\,Miatello}
\address{Facultad de Matem\'atica, Astronom\'\i a y F\'\i sica,
Universidad Nacional de C\'or\-do\-ba, C\'or\-do\-ba~5000, Argentina}
\email{miatello@mate.uncor.edu}

\author{M.I.\,Pacharoni}
\address{Facultad de Matem\'atica, Astronom\'\i a y F\'\i sica,
Universidad Nacional de C\'or\-do\-ba, C\'or\-do\-ba~5000, Argentina}
\email{pacharon@mate.uncor.edu}

\date{\datum}
\maketitle

\section{Introduction} Let $F$ be a totally real number field of
dimension $d$, and let $\Okr_F$ be its ring of integers. If $\qq$ is
an ideal in $\Okr_F$ let $\G=\G_0(\qq)$ denote the congruence
subgroup of Hecke type of the Hilbert modular group. In \cite{BMP} we
have proved a sum formula of Kuznetsov type for $\G$, in which all
weights contribute, and have applied it to give estimates for
averages of Kloosterman sums for $F$. That sum formula has the
following type:
\begin{align*}\int_{Y} \prod_{j=1}^d k_j(\n_j) \, d\sigma_{r,r}(\nu) =
\int_{Y} \, \prod_{j=1}^d k_j(\nu_j)\, d \dt(\nu) + K(\B k)
\end{align*}
 Here, the index $j$ runs over the infinite places of~$F$, the test
functions $k_j$ are even and holomorphic on a strip in~$\CC$ and the
set $Y$ of spectral parameters is a subset of $\CC^d$. The measure $d
\dt$ on~$Y$ has an elementary description (see \eqref{Deltadef}
and~\eqref{Etadef}). The measure $d\s_{r,r}$ is supported on the set
of spectral parameters of automorphic representations, and has
weights that are essentially products of Fourier coefficients of
automorphic forms for~$\G$ (see \eqref{measure}, \eqref{dmeasure},
\eqref{cdef}). The term $K(\B k)$ is a sum of Kloosterman sums,
depending on a Bessel transform $\B k$ of $k=\times_j k_j$.

In this paper we take a proper subset $E\subset\vz{1,\ldots,d}$,
choosing $k_j$ conveniently for $j\in Q:=\vz{1,\ldots,d}\setminus E$
and leaving $k_j$ free for $j\in E$. This choice will lead to a
partial sum formula involving test functions of product type at the
places in $E$ (Theorem \ref{propv}). To prove this, we apply the sum
formula in \cite{BMP} to suitable test functions depending on a
 parameter, and estimate the contributions of the different terms in
the formula. The main task will be to show that the contributions of
 the so called Kloosterman and Eisenstein terms are of lower order of
 magnitude than that of the delta term (see Sections 4 and 5). For the
 estimation of the Eisenstein term, we will need to give an estimate
on the vertical line $\re\n=0$, for Fourier coefficients of
Eisenstein series at each cusp of $\G_0(\qq)$, making explicit the
dependence on the order~$r$ of the Fourier term. In the estimate of
the Fourier coefficient, we will use a logarithmic lower bound for
ray class $L$-functions on the critical line. This will be obtained
by an argument similar to one given by Landau in the case of the
Dedekind zeta function (\cite{La}).

The main goal of this paper is the application of the partial sum
formula to derive density results for automorphic representations. We
now describe some consequences of the main result,
Theorem~\ref{mainthm}. We consider first the set ${\mathcal D_X}$ of
those automorphic representations $\irr = \times_1^d \irr_j$ of
$G\simeq SL_2(\RR)^d$ that have a prescribed discrete series
eigenvalue at each place $j \ne l$, and also a discrete series
eigenvalue at $j=l$, with $-X\leq -\ld_{\irr_l}\leq 0$. We shall see
that this set has positive density, that is, the quantity
$X^{-1}\mu_r({\mathcal D_X})$ tends to a positive limit as $X\mapsto
\infty$, where $\mu_r$ denotes a suitable measure in the unitary dual
of $G$, depending on the Fourier coefficients of order~$r$ of
automorphic representations (see Proposition \ref{dseries}).

We also consider, in Proposition~\ref{pseries}, a set ${\mathcal P}_X$
of automorphic representations having eigenvalue parameters that lie
in fixed intervals $I_j\subset \RR ^{\ge 0}$ at all places $j\ne l$,
and such that $\frac 14<\ld_{\irr_l}<X$. In the case when $I_j \cap
(\frac 14, \infty) \ne \emptyset$ for each $j$, i.e. if all
components are of principal series type, we will show that
$X^{-1}\mu_r({\mathcal P}_X)$ tends to a positive constant. This
constant depends on the measure of the intervals and on $F$, but not
 on $r$, nor on $\qq$. Something similar happens with the limit
constant in the case of the set ${\mathcal D_X}$. These results imply
that there are infinitely many automorphic representations such that
all components are either of principal series type or of discrete
series type and have a non-zero Fourier coefficient of order~$r$.

On the other hand, we will see in Proposition~\ref{cseries}, that if
at least one of the intervals $I_j$ is contained in $[0,1/4]$, then
the limit $X^{-1}\mu_r({\mathcal P}_X)$ is zero. This implies that
the automorphic representations which are of complementary series
type at least at one place, are rare, i.e. they have density zero
with respect to the measure $\mu_r$. We note that this is expected,
given the Ramanujan-Petersson conjecture which predicts no components
of complementary series type in automorphic representations.

As another main application, if we consider the set of all automorphic
representations ${\mathcal A}_X$ with eigenvalue components
$\ld_{\irr_j}$ satisfying the condition $ \sum_{j=1}^d
|\ld_{\irr_j}|< X$, then we have that $X^{-d}\, \mu_r({\mathcal
A}_X)$ tends to a positive limit as $X\mapsto \infty$ (see Corollary
\ref{Eempty}). This result can be seen as similar to a Weyl law,
weighted by Fourier coefficients of automorphic representations.

A Weyl law for spherical automorphic representations for congruence
subgroups of the Hilbert modular group was proved by Efrat
(\cite{Ef}), by using the Selberg trace formula on $\G\backslash G/K$
(\cite{Efmem}). This result implies the existence of infinitely many
$K$-spherical automorphic representations in this context. This
 existence result follows from Theorem~\ref{mainthm}. Efrat counts
vectors of eigenvalues by their $l^2$-norm, whereas our result uses
the $l^1$-norm. Moreover, our distribution results have squares of
Fourier coefficients as weights.

Results related to those in this paper, in the case when $d=1$, were
obtained by Bruggeman (\cite{Br78}, \S4) and by Deshouilliers-Iwaniec
(see~\cite{DI}, Theorem 2). In the case of the Lie group ${\rm
SU}(2,1)$ and the trivial $K$-type, a similar result was given by
Reznikov in \cite{Re93}.
\smallskip

The authors wish to thank F.~Shahidi for very useful conversations on
the most recent best bounds on exceptional eigenvalues.

\section{Preliminaries}\label{prelims} As in~\cite{BMP}, let $F$ be a
totally real number field, and let $\Okr$ be its ring of integers. We
consider the algebraic group $\GG=R_{F/\QQ}(\SL_2)$ over~$\QQ$
obtained by restriction of scalars applied to $\SL_2$ over~$F$.

Let $\s_1,\ldots, \s_d$ be the embeddings $F\rightarrow\RR$. We have
\begin{equation} \label{Gdef}
G\isdef \GG_\RR \cong \SL_2(\RR)^d, \qquad \GG_\QQ \cong \vzm{
\left(x^{\s_1},\ldots, x^{\s_d} \right) }{x\in \SL_2(F) },
\end{equation}
$G$ contains $K\isdef \prod_{j=1}^d \SO_2(\RR)$ as a maximal compact
subgroup.

The image of $\SL_2(\Okr)\subset \SL_2(F)$ corresponds to $\GG_\ZZ$.
This is a discrete subgroup of $\GG_\RR$ with finite covolume. It is
called the {\sl Hilbert modular group}, see \cite{Fr}, \S3. We choose
a non-zero ideal~$\qq$ in~$\Okr$ and form the congruence subgroup of
Hecke type $\G=\G_0(\qq) = \vzm{\matc abcd\in \SL_2(\Okr)}{c\in
\qq}$, which has finite index in~$\GG_\ZZ$.

In this paper, we are concerned with functions on $\G\backslash G$. We
restrict ourselves to \emph{even functions:} $f(-g)=f(g)$. By
$L^2(\G\backslash G)^+$ we mean the Hilbert space of (classes of)
even functions that are left invariant under~$\G$, and square
integrable on $\G\backslash G$ for the measure induced by the Haar
measure. This Hilbert space contains the closed subspace
$L^2_c(\G\backslash G)^+$ generated by integrals of Eisenstein
series. The orthogonal complement $L^2_d(\G\backslash G)^+$ of
$L^2_c(\G\backslash G)^+$ is the closure of $\sum_{\irr} V_\irr$,
where $V_\irr$ runs through an orthogonal family of closed
irreducible subspaces for the $G$-action in~$L^2(\G\backslash G)$ by
right translation.

\rmk{Irreducible unitary representations}Each representation $\irr$
has the form $\irr= \otimes_j \irr_j$, with $\irr_j$ an even unitary
irreducible representation of~$\SL_2(\RR)$. Table~\ref{irrrepr} lists
the possible isomorphism classes for each~$\irr_j$. For each~$\irr$
we define a spectral parameter $\n_\irr =
\left(\n_{\irr,1},\ldots,\n_{\irr,d}\right)$, with $\n_{\irr,j}$ as
in the last column of the table. There are Casimir operators $C_j$
acting on each coordinate for $1\le j \le d$. The
eigenvalue~$\ld_\irr\in \RR^d$ is given by $\ld_{\irr,j} =
\frac14-\n_{\irr,j}^2$. We note that if $\irr_j$ lies in the
complementary series, $\ld_{\irr,j} \in (0,\frac 14)$ and if $\irr_j$
is isomorphic to a discrete series representation $D^{\pm}_b$, $b \in
2\ZZ$, $b\ge 2$, then $\ld_{\irr,j}=\frac b2 (1- \frac b2)\in \ZZ_{
\le 0}$. If $Q$ is an arbitrary subset of $\vz{1,\ldots,d}$, then we
shall denote by ${\| \ld_{\irr,Q} \|}_1=\sum_{j\in Q}
|\ld_{\irr,j}|$, the 1-norm of the projection $\ld_{\irr,Q}$, of
$\ld_{\irr}$ onto the subspace of $\RR^d$ corresponding to $Q$.

\begin{table}\small
\[ \widearray{|cc|c|c|c|}\hline
\multicolumn2{|c|}{\text{notation}}&\text{name}&\text{weights}&\n\\
\hline
1&&\text{trivial representation}&0&\frac12\\
H(s)&s\in i[0,\infty)&\text{unitary principal series}&q\in 2\ZZ& \frac
s2 \\
H(s)&s \in(0,1)&\text{complementary series}&q\in2\ZZ& \frac s2 \\
D^+_b&b\geq2, \,b\in 2\ZZ&\text{holomorphic discrete series}&q\geq
b,\, q\in 2\ZZ&\frac{b-1}2\\
D^-_b&b\geq2,\, b\in 2\ZZ&\text{antiholomorphic discrete series}&q\leq
-b,\, q\in 2\ZZ& \frac{b-1}2
\\
\hline
\end{array}
\]
\caption[i]{{\em Irreducible unitary even representations of the Lie
group $\SL_2(\RR)$.} All characters of~$SO_2(\RR)$ occur at most
once; the characters that occur are listed under {\sl weights}. The
last column gives a spectral parameter $\n$, with $\re\n\geq0$, such
that $\ld(\n)=\frac14-\n^2$ is the eigenvalue of the Casimir
operator. See~\cite{La75}, Chap.~VI, \S6.}\label{irrrepr}
\vspace{.1cm}
\end{table}

The constant functions give rise to $\irr = \one \isdef \otimes_j 1$.
It occurs with multiplicity one. If $V_\irr$ does not consist of the
constant functions, then $\irr_j\neq 1$ for all~$j$.

If $\n_j\in (0,\frac12)$ for some~$j=1,\ldots,d$, then we call
$\ld(\n):=\frac14-\n^2 \isdef \left( \frac14-\n_j^2\right)_j$ an {\em
exceptional eigenvalue}. We call such a $\n_j$ an {\em exceptional
coordinate}. If $d=1$, it is known that only finitely many
exceptional eigenvalues can occur for a given~$\G$. For $d>1$ such a
result has not been proved. In principle, there might be infinitely
many exceptional eigenvalues, since one coordinate can stay small,
while others tend to $\infty$. On the other hand, the
Ramanujan-Petersson conjecture predicts that there are none, that is,
$\nu_j \in i\RR\cup \left(\frac12+\ZZ_{\geq0}\right)$, for all $j$.
We note that the results of Efrat, \cite{Ef}, imply that there exist
only finitely many eigenvalues such that {\it all} coordinates are
exceptional.

In the case $F=\QQ$, Selberg showed that $\n_j
\not\in\left(\frac14,\frac12\right]$ for such exceptional
coordinates, that is $\ld(\n_j)\ge\frac{3}{16}$ (\cite{Se65}). This
estimate has been improved and extended to arbitrary number fields by
several authors. In the general case, the best bound in the
literature (see~\cite{LRS}) is $\nu_j \le \frac 15$, that is,
$\ld(\n_j)\ge \frac{21}{100}$, for $1\le j\le d$. This bound has
recently been improved to $\ld(\n_j)\ge 0.22837$; see \cite{KS1} and
\cite{KS2}.

\rmk{Automorphic forms and Fourier coefficients} The elements
of~$V_\irr$ that transform on the right according to a character of
the maximal compact subgroup~$K$ are square integrable automorphic
forms. The Fourier coefficients of the automorphic forms are
essentially independent of the actual choice of the automorphic forms
in~$V_\irr$, and are determined by the number
$c^r(\irr)=c_\infty^r(\irr)$ in Equation~(17) of~\cite{BMP}. As we
consider only the cusp~$\infty$, we omit the parameter~$\k$ from the
notation. The number $r\in \Okr'=\vzm{ x\in F}{\Tr_{F/\QQ}(xy)\in\ZZ
\text{ for all }y\in\Okr}$ determines the order of the Fourier
coefficient.

We take $r\neq0$. That implies that $\irr\neq\one$ if $c^r(\irr)\neq0$
for some~$r\neq0$. If $q=(q_1,\cdots,q_d)\in2\ZZ$ is a weight
occurring in~$\irr$, then there is an $\ps\in \irr$, normalized as in
\cite{BMP}, (15), with the following Fourier term of order~$r$ in the
point $g = \left( \matc{\sqrt{y_j}}{x_j/\sqrt{y_j}}0{1/\sqrt{y_j}}
\matr{\cos\th_j}{\sin\th_j}{-\sin\th_j}{\cos\th_j} \right)\in G$:
\begin{align}\label{cdef}
&\frac{c^r(\irr) }{\vol((\G\cap N)\backslash N)}
\\
\nonumber
&\quad\hbox{} \cdot \prod_{j=1}^d e^{2\p ir_jx_j} \frac{ (-1)^{q_j/2}
(2\p|r_j|)^{-1/2} }{ \G\left(\frac12+\n_j + \frac12 q_j\sign(r_j)
\right) } W_{\sign(r_j)q_j/2,\n_j}(4\p|r_j|y_j) e^{iq_j\th_j} ,
\end{align}
where $N= \left\{ \matc 1x01 \right\}\subset G$, and
$W_{\cdot\,\cdot}$ is the exponentially decreasing Whittaker
function.

\begin{defn} \emph{Test functions. }\label{Lkrdef}
Fix $\tau\in \left(\frac12,\frac34\right)$. Let $\Lkr$ be the space of
even holomorphic functions on the set
\[
\vzm{\n\in\CC}{|\re\n|\leq\tau}\cup\textstyle{\left(\frac12+\ZZ\right)}\]
that satisfy the conditions
\begin{align*}
&k(\n) \ll \left(1+|\im\n|\right)^{-a} \text{ for } |\re\n|\leq\tau
\text{ for some } a>2, \text{ and } \\
&\sum_{b\in 2\ZZ,\, b\geq2} {\textstyle\frac{b-1}2}
\left|k\left({\textstyle\frac{b-1}2}\right)\right|<\infty.
\end{align*}
\end{defn}

For each such test function, we define the following quantity:
\begin{equation} \label{Etadef}
\Eta (k) \isdef \tfrac i2\int_{\re\n=0} k(\n)\, \n\,\tan\p\n\, d\n +
\sum_{b\geq2,\, b\in2\ZZ} {\textstyle\frac{b-1}2
k\left(\frac{b-1}2\right)}.
\end{equation}

\noindent We set
\begin{equation}
\Akr_d=\vzm{\textstyle{\frac{b}2\left(1-\frac b2 \right)}}{b \text{
even }, b\ge2} \subset\ZZ_{\le 0}, \qquad \Ykr:= (0,\infty)\cup
\Akr_d \subset \RR.
\end{equation}
The eigencoordinates $\ld_{\irr,j}$ are elements of $\Ykr$. The test
functions $k\in \Lkr$ give rise to functions of the form
\begin{equation}\label{g}
g(\tfrac 14-\nu^2)=k(\n).
\end{equation}
The domain of these functions consists of the $\frac14-\n^2\in\CC$
with $|\re\n|\leq\tau$, together with the numbers $\frac
b2\left(1-\frac b2\right)$, $b\in\ZZ$. We denote by $\tilde\Lkr$ the
class of functions~$g$ obtained in this way. If $g \in \tilde\Lkr$ we
set
\begin{equation}\label{Etatil}
\tilde\Eta(g)=\tfrac 12\int_{\frac14}^\infty g(y) \tanh (\pi \sqrt{
y-1/4})\, dy +\sum_{y\in \Akr_d} \textstyle\sqrt{1/4-y}\, g(y).
\end{equation}
The functional $\tilde\Eta$ on $\tilde\Lkr$ is described by a measure
$d\eta$ on $\Ykr$, which is given by $d\eta(y) =\frac 12
\tanh\left(\pi\sqrt{y-1/4}\right)\, dy$ on $(1/4,\infty)$, and by the
sum of $\sqrt{1/4-\nobreak y}$ times the delta measure at the points
$y\in \Akr_d$.

\section{Statement of main results}\label{results}
Theorem~\ref{propv} and Theorem~\ref{mainthm} state the main results
of this paper. We give some consequences in
Propositions~\ref{cseries}--\ref{pseries} and Remark~\ref{d1case}. We
 keep the notation from Section~\ref{prelims}.

All statements in this section depend on the partition of
$\vz{1,\ldots,d}$ into three disjoint subsets $E$, $Q_+$ and $Q_-$,
with $Q:=Q_+\cup Q_-\neq\emptyset$.

 To each such partition we attach a set of irreducible representations
 of $G$:
\begin{equation} \label{partition}
\Rkr(E,Q_+,Q_-)= \left\{ \irr\neq\one : \ld_{\irr,j}\geq 0 \text{ if
}j\in Q_+, \ld_{\irr,j}<0, \text{ if }j\in Q_-\right\}.
\end{equation}
We shall often write $\Rkr=\Rkr(E,Q_+,Q_-)$. Our main results will
follow from the next theorem:

\begin{thm}\label{propv}
Let $r\in\Okr'\setminus\vz0$. Choose a partition $E$, $Q_+$, $Q_-$ of
the set $\vz{1,\ldots,d}$ with $Q =Q_+ \cup Q_- \ne \emptyset$ and
let $\Rkr$ be the corresponding set as in (\ref{partition}). If
$g=\times_{j\in E}\, g_j$, with $g_j\in \tilde\Lkr$ for each $j\in
E$, then the series
\begin{equation} \label{Zeta} Z_s(g) \isdef \sum_{\irr\in \Rkr}
\left|c^r(\irr)\right|^2 e^{-s{\| \ld_{\irr,Q} \|}_1} \prod_{j\in E}
g_j(\ld_{\irr,j})
\end{equation}
converges absolutely for each $s>0$, and
\begin{equation}
\label{limZeta} \lim_{s\downarrow 0} s^{d-|E|} Z_s(g)
={\textstyle\frac{ 2^{1+|E|}}{(2\pi)^d}}\sqrt{|D_F|}\, \prod_{j\in E}
\tilde\Eta(g_j)
\end{equation}
Here $\tilde\Eta(g)$ is as given in \eqref{Etatil}. If $E=\emptyset$
then the products over $j\in E$ in the right hand side of
\eqref{Zeta} and \eqref{limZeta} are interpreted as $1$.
\end{thm}
\medskip

The theorem will be proved in Section~\ref{discrterm}, as an
application of the sum formula of Kuznetsov type given in \cite{BMP}.
In this section we shall use Theorem \ref{propv} to obtain density
results for automorphic representations. We first prove the following
proposition:
\begin{prop}\label{thmX}Let $r\in\Okr'\setminus\vz0$. Let $E$, $Q_\pm$
and $\Rkr$ be as above. If $g_j\in \tilde \Lkr$ for $j\in E$, then
\begin{equation}
\begin{split}\label{mainth1}
\lim_{X\rightarrow\infty} X^{|E|-d}& \sum_{\substack{\irr\in\Rkr \\
 {\| \ld_{\irr,Q} \|}_1\leq X}} \, \left| c^r(\irr)\right|^2
\prod_{j\in E} \,g_j(\ld_{\irr,j})\\
&= {\textstyle\frac{ 2^{1+|E|}\,\sqrt{|D_F|}}{(d-|E|)!\, (2\pi)^d}}\,
\prod_{j\in E} \tilde \Eta(g_j),
\end{split}
\end{equation}
with $\tilde\Eta(g_j)$ as given in \eqref{Etatil}. If $E=\emptyset$
then the products over $j \in E$ in \eqref{mainth1} equal $1$.
\end{prop}

\begin{proof}
For each $X>0$, and each $g=\times_{j\in E} \,g_j \in \tilde
 \Lkr^{|E|}$, we define
\begin{equation}\label{muXg}
\mu_g(X)= \sum_{\substack{\irr\in \Rkr
 \\ {\| \ld_{\irr,Q} \|}_1 \leq X}} \left|c^r(\irr)\right|^2\,
 \prod_{j\in E} \,g_j(\ld_{\irr,j}).
\end{equation}
 The absolute convergence follows from Theorem~\ref{propv}.

If $Z_s(g)$ is as in~\eqref{Zeta}, we have
\begin{equation}
 Z_s(g)=\sum_{\irr\in \Rkr } \left|c^r(\irr)\right|^2 e^{-s{\|
\ld_{\irr,Q} \|}_1}\,\prod_{j\in E} g_j(\ld_{\irr,j}) =
\int_{X=0}^{\infty} e^{-sX }\, d\mu_g(X).
\end{equation}
By applying Theorem \ref{propv}, we find that
\begin{equation}
\lim_{s\rightarrow 0} s^{d-|E|} \int_{X=0}^{\infty} e^{-sX }\,
d\mu_g(X) = {\textstyle\frac{ 2^{1+|E|}}{(2\pi)^d}}\sqrt{|D_F|}
\,\prod_{j\in E} \tilde\Eta(g_j).
\end{equation}
\smallskip

In the case when the functions $g_j$ are non-negative on $\Ykr=
(0,\infty) \cup \Akr_d$, then the function $X\mapsto \mu_g(X) $ is
non-decreasing and we may apply a Tauberian theorem (see Theorem~4.3
in Chap.~V of~\cite{Wi}) to obtain:
\begin{equation}\label{limmuX}
\lim_{X \rightarrow \infty} X^{|E|-d} \mu_g(X) = {\textstyle\frac{
2^{1+|E|}\,\sqrt{|D_F|}}{(d-|E|)!\,(2\pi)^d}}\, \prod_{j\in E} \tilde
\Eta(g_j)
\end{equation}

For general $g_j \in \tilde \Lkr$, the bounds in
Definition~\ref{Lkrdef} allow us to construct $\tilde g_j\in \tilde
\Lkr$ such that $\tilde g_j\ge 0$ on $\Ykr$ and $|g_j|\le \tilde
g_j$.

Let $\g(X)$ be equal to $\re \mu_g(X)$ or $\im \mu_g(X)$. For $X_1 >X$
we have:
\begin{equation*}
\g(X_1) -\g(X)= \sum_{\substack{\irr\in\Rkr
\\ X<{\| \ld_{\irr,Q} \|}_1\leq X_1}} \left| c^r(\irr)\right|^2 \, (
\re \text{ or } \im) \prod_{j\in E} \,g_j(\ld_{\irr,j})
\end{equation*}
hence $|\g(X_1) -\g(X)| \le \mu_{\tilde g} (X_1) -\mu_{\tilde g}(X).$
So $\bt (X)\isdef \g (X) + \mu_{\tilde g}(X)$ is non-decreasing and
we have
\begin{align*}
\lim_{X\mapsto \infty} X^{|E|-d}&\bt(X) \\
=& {\textstyle\frac{ 2^{1+|E|}}{(d-|E|)!(2\pi)^d}}\sqrt{|D_F|}\,
\left( (\re \text{ or } \im )\prod_{j\in E}\tilde\Eta(g_j)+
\prod_{j\in E}\tilde\Eta(\tilde g_j)\right).
\end{align*}
Since \eqref{limmuX} holds for $\tilde g = \prod_{j\in E}\tilde g_j$,
this implies that \eqref{limmuX} holds also for arbitrary $g_j\in
\tilde \Lkr$, and the proposition follows.
\end{proof}

We now state the main result in this paper. In the proof we will
extend the validity of \eqref{mainth1} to a larger class of functions
$g$ than those considered so far.

\begin{thm}\label{mainthm} Let $r\in\Okr'\setminus\vz0$. Let
$\vz{1,\ldots,d}$ be the disjoint union of the subsets $E$, $Q_+$ and
$Q_-$, with $E\neq\vz{1,\ldots,d}$. Let $H$ be a hypercube
$H=\prod_{j\in E} [a_j,b_j] \subset \RR^{|E|}$, such that $a_j,b_j
\not\in \vzm{ \frac b2(1-\frac b2)}{b\geq2,\, b\text{ even}}$. Then,
if $\Rkr=\Rkr(E,Q_+,Q_-)$ is as in~\eqref{partition}, we have
\begin{equation}\label{dens}
\begin{split}
\lim_{X\rightarrow\infty} &X^{|E|-d} \sum_{ \substack{\irr\in\Rkr ,\,
a_j\leq \ld_{\irr,j}\leq b_j,\, j\in E
\\ {\| \ld_{\irr,Q} \|}_1\leq X}} \, \left| c^r(\irr)\right|^2
\\
&= {\textstyle\frac{ 2 \,\sqrt{|D_F|}}{(d-|E|)!\, (2\pi)^d}}\,
\prod_{j\in E} \Bigg( \int_{[a_j,b_j]\cap [1/4,\infty)}
\tanh\left(\pi \sqrt{y-1/4}\right)\, dy
\\
&\quad \qquad \qquad \qquad \hbox{} + \sum_{\substack{b\geq 2,\,
b\,{\rm even},\\ a_j< \frac b2(1-\frac b2)<b_j} } (b-1) \Bigg).
\end{split}
\end{equation}
\end{thm}
\rmk{Remark}Note that the factors in the product over $j\in E$ in the
right hand side of \eqref{dens} are twice the volume of $[a_j,b_j]$
for the measure $d\eta$ discussed after~\eqref{Etatil}.

\begin{proof} As a first step in the proof we will show that
\eqref{mainth1} is also valid for any function $k(\ld)= \times_{j \in
E} k_j (\ld)$ with $k_j (\ld)$ an arbitrary continuous, compactly
supported and real valued function on $\tilde \Ykr$. For this
purpose, we carry out an approximation argument based on the fact
that for each $\e>0$ there exists $h_\e=\times_{j\in E}\, h_j \in
\tilde\Lkr^{|E|}$ such that
\begin{equation}\label{ghb}
\left| k(\ld) - h_\e(\ld ) \right|\leq \e b(\ld)\quad\text{ for all }
\ld \in \Ykr^{|E|},
\end{equation}
where $b=\times_{j\in E}\, b_j$ is also an element of
$\tilde\Lkr^{|E|}$, not depending on~$\e$ with $b_j>0$, for each $j$.

We first show how \eqref{ghb} leads to the assertion. After that, we
construct $b$ and $h_\e$ satisfying \eqref{ghb}. In the remainder of
the proof we shall write $\mu_X(g):=\mu_g(X)$ to stress the
dependence on the test function $g$. We will use the fact that
$\mu_X$ defines a positive measure on $\Ykr^{|E|}$ for each fixed
$X$. Also, we denote by $\m$ the non-negative measure on $\Ykr^{|E|}$
in the right hand side of~\eqref{mainth1}, given by
\begin{equation*}
\int_{\Ykr} f(y)d\m(y) ={\textstyle\frac{
2^{1+|E|}\,\sqrt{|D_F|}}{(d-|E|)!\,(2\pi)^d}}\,
\int_{\Ykr^{|E|}}f((y_j)_{j\in E}) \prod_{j\in E}{d\eta(y_j)}.
\end{equation*}
We have
\[ h_\e(\ld) - \e b(\ld) \leq k(\ld) \leq h_\e(\ld) + \e b(\ld)\]
for $\ld\in \Ykr^{|E|}$. Hence
\begin{alignat}2 \label{muhX}
\m_X(h_\e)-\e \m_X(b) &\leq& \m_X(k) &\leq \m_X(h_\e) + \e\m_X(b),\\
\label{muh}
\m(h_\e)-\e \m(b) &\leq& \m(k) &\leq \m(h_\e) + \e\m(b).
\end{alignat}

The inequality~\eqref{muh} shows that
\[ \left|\m(k) - \m(h_\e)\right|\leq \e \m(b) \]
for each $\e>0$. Hence $\lim_{\e\downarrow0}\m(h_\e) = \m(k)$.

Proposition~\ref{thmX} shows that $\lim_{X\rightarrow\infty} X^{|E|-d}
\m_X (h_\e) = \m(h_\e)$ for each $\e>0$, and
$\lim_{X\rightarrow\infty} X^{|E|-d} \m_X (b) = \m(b)$. We want to
prove the same for~$k$.

We derive from~\eqref{muhX} that
\[0\leq \limsup_{X\rightarrow\infty} X^{|E|-d} \m_X(k)
-\liminf_{X\rightarrow\infty} X^{|E|-d} \m_X(k) \leq 2\e\m(b).\]
Hence $\lim_{X\rightarrow\infty} X^{|E|-d} \m_X(k) $ exists; by taking
limits in~\eqref{muhX} we obtain:
\[ \m(h_\e) - \e\m(b) \leq \lim_{X\rightarrow\infty} X^{|E|-d} \m_X(k)
\leq \m(h_\e) + \e\m(b).\]
The desired equality follows by taking the limit as $\e\downarrow0$.
\medskip

Now we turn to the construction of $h_\e$ and $b$. First we focus our
attention on one place $j\in E$. For positive $u$, we put
\begin{align} \label{hjdef}
 h_j(\ld) &\isdef \sqrt{\frac{u}\p} \int_{-\infty}^\infty e^{-uy^2}
 k_j(\ld-y)\, dy\qquad (\ld\in\RR).
\end{align}
The dependence of $h_j$ on~$u$ is not visible in the notation. This
expression shows that $\n\mapsto h_j(1/4-\nobreak\n^2)$ is
holomorphic on~$\CC$, and has an exponential decay on the strip $|\re
\n| <\tau$ and along the real axis. This implies that $h_j\in
\tilde\Lkr$. We note that for each $j$
\begin{equation}\label{hsmtk}
{\| h_j \|}_\infty \le {\| k_j \|}_\infty
\end{equation}
The function $h_j$ gives a holomorphic approximation of~$k_j$. By
taking $u$ sufficiently large, we obtain, for prescribed $\e_j>0$,
\begin{equation}\label{hmksmte}
\|h_j-\nobreak k_j\|_\infty <\e_j.
\end{equation}
Let us take $A_j> 1$ so that $\supp (k_j) \subset [-A_j +1 , A_j -1]$.
If $|\ld|>A_j$ then we have:
\begin{align*}
|h_j(\ld) | &\leq 2(A_j-1) \|k_j\|_\infty \sqrt{u/\p}
\;e^{-u(|\ld|+1-A_j)^2} \\
&\leq 2(A_j-1) \|k_j\|_\infty \sqrt{u/\p}\; e^{-u(|\ld|+1-A_j)}
\end{align*}
Further enlarging~$u$, if necessary, we get $e^{-u} \sqrt{u/\p}
\;2(A_j -1)\|k_j\|_\infty \le \e_j $, and
\begin{equation}\label{hsmteb}
|h_j(\ld)|\leq \e_j b_j(\ld),\;\; \text{ for } |\ld|\geq A_j,
\end{equation}
 where
\[ b_j(\ld) \isdef \begin{cases}
\left( \frac{1+\ld}{1+A_j}\right)^{-2} &\text{ for }\ld\geq
-\frac12,\\
\left( \frac{1-\ld}{1+A_j}\right)^{-2} &\text{ for }\ld <-\frac12.
\end{cases}
\]
We extend $b_j$ as a holomorphic function on a neighborhood of
$\;[0,\infty)$, such that $\n\mapsto b_j(1/4-\nobreak\n^2)$ is
holomorphic on $|\re\n|\leq \tau$, with $\tau$ slightly larger
than~$\frac12$. The estimate $b_j(1/4-\nobreak
\n^2)=\oh\left((1+|\n|)^{-4}\right)$ is sufficient to conclude that
$b_j\in\tilde\Lkr$.

Note that the $h_j$ depend on $\e_j>0$ but the $b_j$ do not. If $| \ld
|\le A_j$, then we have that $b_j (\ld)\ge 1$ and we conclude that
\begin{equation}\label{hsmtbtb}
|h_j(\ld)| \leq\; \bt_j b_j(\ld)
\end{equation}
where $\bt_j= 2(A_j-1) \|k_j\|_\infty \sqrt{u/\p}$.

We now take $h_\e=\times_{j\in E}\, h_j$ and $b=\times_{j\in E} \,
b_j$, and assume that all $\e_j\in (0,1)$. For a given $\ld\in
\Ykr^{|E|}$, let $F=\vzm{j\in E}{|\ld_j|\leq A_j}$. If $F=E$, then we
have by \eqref{hsmtk}, \eqref{hmksmte}, and the fact that
$b_j(\ld_j)\geq1$ if $|\ld_j|\leq A_j$:
\begin{align*}
\left|k(\ld)-h_\e(\ld)\right| &\leq \sum_{j\in E} \left| k_j(\ld_j) -
h_j(\ld_j)\right| \prod_{i\in E,\, i\neq j} \|k_i\|_\infty
\\
&\leq \Big( \sum_{j\in E}\e_j\Big) \prod_{i\in E} \max\left( 1,
\|k_i\|_\infty\right) \; b(\ld).
\end{align*}
The product over $i\neq j$ is taken outside the sum over~$j$ by
estimating it by \[\prod_{i\in E} \max\left( 1,\|k_i\|_\infty\right)
\; b(\ld).\]
For $F\neq E$, we use \eqref{hsmtbtb} and \eqref{hsmteb},
to obtain:
\[ \left|k(\ld)-h_\e(\ld)\right|=|h_\e(\ld)| \leq b(\ld) \prod_{j\in
F}\bt_j \prod_{i\in E\setminus F} \e_i .\]
Thus we see that we can adjust the $\e_i$ in such a way that
condition~\eqref{ghb} is satisfied.
\medskip

As the final step in the proof of the theorem we extend
\eqref{mainth1} to the characteristic function of a hypercube
$\prod_{j\in E} [a_j,b_j] \subset \RR^{|E|}$.

 Let us denote by $\chi_j$ the characteristic function of $[a_j,
b_j]$. Now, for any $\e >0$, it is easy to construct functions $u_j,
U_j \in C_c(\RR)$ such that $ 0\leq u_j \le \chi_j \le U_j$, $
\parallel
u_j\parallel_\infty\le 1$, $ \parallel U_j\parallel_\infty\le 1$, and
$ \int_{-\infty}^{\infty} (U_j-u_j) \,d\eta \le \e$. The measure
$d\eta$ has point masses at the elements of $\Akr_d$. So we need the
assumption $a_j,b_j\not\in \Akr_d$ to attain the last inequality for
all $\e>0$.

If we let $u=\times_{j \in E}\,u_j$ and $U=\times_{j \in E}\,U_j$, we
have
\begin{equation}\label{uU}
\m(U)-\m(u) \le C_1\,|E| \e \qquad \text{ and } \qquad u\le \chi \le
U.
\end{equation}
with $C_1$ a constant depending on the supports of the $U_j$ and on
the measure~$\m$.

Now since $\m_X(u)\le \m_X(\chi) \le \m_X(U)$ and \eqref{mainth1} is
valid for $u$ and $U$ we see that
\begin{equation}\label{mulim}
\m(u)\le \liminf_{X\rightarrow \infty}X^{|E|-d}\m_X(\chi) \le
\limsup_{X\rightarrow \infty}X^{|E|-d}\m_X(\chi) \le \m(U)
\end{equation}
The existence of $\lim_{X\rightarrow \infty}X^{|E|-d} \m_X(\chi)$
follows from \eqref{uU} and \eqref{mulim}. Since $\m(u)\le \m(\chi)
\le \m(U)$, we have that
$$\lim_{X\rightarrow \infty}X^{|E|-d}\m_X(\chi)=\m(\chi).$$
This is the statement in the theorem.
\end{proof}

If we let $E=\emptyset$ in Theorem \ref{mainthm}, we get the following
result that can be seen as a Weyl-type law, weighted by Fourier
coefficients.

\begin{cor} \label{Eempty} Let $r\in \Okr'\setminus\vz0$ and let
$\vz{1,\ldots,d} = Q_+\sqcup Q_-$. Then we have:
\begin{align}
\lim_{X\rightarrow \infty} X^{-d} \sum_{\substack{ \|\ld _{\irr} \|_1
<X\\
\ld_{\irr,j}>0,\, j\in Q_+\\
\ld_{\irr,j}<0,\, j\in Q_- }} \left|c^r(\irr)\right|^2 =
{\textstyle\frac{ 2\,\sqrt{|D_F|}}{d!\, (2\pi)^d}}.
\end{align}
\end{cor}

The possibility to prescribe $Q_+$ and $Q_-$ allows us to count
representations having discrete series type factors at some places,
and factors of principal or complementary series type at the other
places. The following result ignores this distinction.

\begin{cor}\label{corgen}
Let $r\in \Okr'\setminus\vz0$. Let $E$ be a proper subset of
$\vz{1,\ldots,d}$, and put $Q=\vz{1,\ldots,d} \setminus E$. Take
$[a_j,b_j]$, for $j\in E$, as in the theorem. Then
\begin{equation}\label{corlim}
\begin{split}
\lim_{X\rightarrow\infty} & X^{|E|-d} \sum_{ \substack{ a_j\leq
\ld_{\irr,j}\leq b_j,\, j\in E
               \\{\| \ld_{\irr,Q} \|}_1 \leq X } }
\left|c^r(\irr)\right|^2\\
&={\textstyle\frac{ 2\,\sqrt{|D_F|}}{(d-|E|)!\, \pi^d 2^{|E|}}}\,
\prod_{j\in E} \left( \int_{[a_j,b_j] \cap [1/4,\infty)}
 \tanh\left(\pi \sqrt{y-1/4}\right)\, dy \right.
\\
&\quad \qquad \qquad \qquad \left. \hbox{} + \sum_{\substack{ b\geq
2,\, b\,{\rm even}\\ a_j< \frac b2(1-\frac b2)<b_j}} (b-1)
\right.{\Bigg )}.
\end{split}
\end{equation}
\end{cor}
\begin{proof}This is obtained from the theorem by adding the
contributions, over all possible choices of $Q_+\sqcup
Q_-=\vz{1,\ldots,d} \setminus E$.
\end{proof}\medskip

We now show that Theorem~\ref{mainthm} can be used to derive density
results for automorphic representations subject to restrictions at
some places.

The first result (Proposition~\ref{cseries}) confirms that the
representations of complementary series type are rare. We fix one
place $l\in \vz{1,\ldots,d}$, and apply Corollary~\ref{corgen} with
$E=\vz l$, and $[a_l,b_l]\subset \left( 0,\frac14\right]$. We catch
all complementary series eigenvalues at the place~$l$ if we take
$0<a_l\leq \frac{21}{100}$, since, according to~\cite{LRS}, all
complementary series factors correspond to eigenvalues in
$\left[\frac{21}{100},\frac14\right)$ (see also \cite{KS1} and
\cite{KS2}). The right hand side in~\eqref{corlim} vanishes for this
choice, and leads to case~(i) in the following proposition.

To obtain part~(ii), we take $E=\vz{1,\ldots,d}\setminus\vz l$, and
fix $k\neq l$. We take $[a_k,b_k]\subset\left(0,\frac14\right]$, and
for $j\neq k$, we let $[a_j,b_j]$ be any interval in~$\RR$, as in the
theorem.

\begin{prop}\label{cseries}For each $r\in \Okr'\setminus\vz0$, we have
$$(i)\quad \lim_{X\rightarrow \infty} X^{1-d} \sum_{\substack{
\sum_{j\neq l} |\ld_{\irr,j}| \leq X \\ 0\leq \ld_{\irr,l}\leq \frac
14}} \, \left| c^r(\irr)\right|^2 =0$$
$$(ii)\quad \lim_{X\rightarrow \infty} X^{-1}
\sum_{\substack{|\ld_{\irr,l}|\leq X\\ 0 \le |\ld_{\irr,k}|\le \frac
14 \\
 a_j \le |\ld_{\irr,j}|\le b_j,\,\, {j\neq k,l} }} \, \left|
 c^r(\irr)\right|^2 =0$$
\end{prop}
Replacing $\ld_{\irr,l}\leq \frac14$ by $\ld_{\irr,l}<\frac14$, we
obtain a density zero result for exceptional eigenvalues.
\medskip

In the next application, we restrict our attention to discrete series
type eigenvalues, and, moreover, prescribe the eigenvalue at all
places but one. So we choose $E=\vz{1,\ldots,d}\setminus\vz l$ and
 $Q=Q_-= \{ l \}$. For each $j\in E$, we pick $\ld_j\in \Akr_d$, and
choose $[a_j,b_j] $ such that $[a_j,b_j ] \cap \Ykr=
\vz{\ld_j}\subset(a_j,b_j)$, and $b_j<\frac{21}{100}$ if $\ld_j=0$.
Application of Theorem~\ref{mainthm}. gives:
\begin{prop}\label{dseries}Let $r\in \Okr'\setminus\vz0$. Let $1\leq
l\leq d$, and take $\ld_j\in \Akr_d$ for $j\neq l$. Then
$$\lim_{X\rightarrow \infty} X^{-1} \sum_{\substack{-X\le
\ld_{\irr,l}<0 \\ \ld_{\irr,j}=\ld_j,\, j\neq l}}
 \left|c^r(\irr)\right|^2 = {\textstyle\frac{\sqrt{|D_F|}}{\pi^d}}\,
\prod_{j\ne l} \sqrt {1/4-\ld_j}.$$
\end{prop}
This shows that there are infinitely many $\irr$ that have discrete
series type factors at all places, and a prescribed eigenvalue at all
but one place. If we take $r$ totally positive, we restrict the sum
to $\irr$ that are generated by a holomorphic Hilbert modular cusp
form (see Proposition~2.2.3 in~\cite{BMP}). The occurrence of the
Fourier coefficients $c^r(\irr)$ in our result makes it hard to find
a connection to dimension formulas for spaces of holomorphic Hilbert
modular forms like those in Theorem~3.5 of~\cite{Fr}.

The number of places at which we have restricted $\irr$ is reflected
in the exponent of $X$. If $d>2$, the positive density here might
actually correspond to a lower density than the zero density in (i)
of the previous proposition.

\medskip
Finally, we take $E=\vz{1,\ldots,d}\setminus\vz l$, $Q=Q_+ =\{l\}$ and
confine $\ld_{\irr,j}$, $j\neq l$, to a small interval $[a_j,b_j]$ of
principal series eigenvalues. We obtain:
\begin{prop}\label{pseries}Let $r\in \Okr'\setminus\vz0$; let $1\leq l
\leq d$, and take $[a_j,b_j] \subset\left[\frac14,\infty\right]$ for
$j\neq l$. Then
\begin{align*}
\lim_{X\rightarrow \infty} & X^{-1} \sum_{\substack{ 0\le \ld_{\irr,l}
\leq X \\a_j \leq \ld_{\irr,j}\leq b_j,\, j\neq l}}
 \left|c^r(\irr)\right|^2
\\
&= {\textstyle\frac{2^{1-d} \sqrt{|D_F|}}{ \pi^d}}\, \prod_{j\neq l}
\int_{a_j}^{b_j} \tanh\left(\p\sqrt{y-1/4}\right)\, dy\\
&= {\textstyle\frac{2^{1-d} \sqrt{|D_F|}}{ \pi^d}}\, \prod_{j\neq l}
\left( (b_j-a_j)\left(1+
\oh\left(e^{-2\p\sqrt{a_j-1/4}}\right)\right) \right).
\end{align*}
\end{prop}

\smallskip
\begin{remark}\label{d1case} In the case when $d=1$, necessarily
$E=\emptyset$ in Theorem \ref{mainthm} and $Q=Q_+ = \{1\}$ or $Q=Q_-
=\{1\}$. The asymptotic result obtained in each case is expressed in
Corollary \ref{Eempty}. We get, as $X\mapsto +\infty$:
\begin{equation}
\begin{split}\label{d=1}
\sum_{\substack{ 0 \le \ld _{\irr} \le X }} \left|c^r(\irr)\right|^2
\sim {\textstyle\frac{ X}{\pi}},\quad \text{ if } Q_-=\emptyset,\\
\sum_{\substack{ -X \le \ld _{\irr} < 0 }} \left|c^r(\irr)\right|^2
\sim {\textstyle\frac{ X}{\pi}},\quad \text{ if } Q_+=\emptyset.
\end{split}
\end{equation}
We note that in \cite{Br78}, Corollary 4.4 gives a special case of
Theorem~\ref{propv}, implying the first asymptotic formula
in~\eqref{d=1}. Also, in \cite{DI}, Theorem 2, a result that applies
 to Fourier coefficients of general automorphic forms is given in the
 form of an upper bound.

In comparing the result for $d=1$ with those in \cite{Br78} and
\cite{DI}, it is useful to note that a normalized with respect to the
usual measure $y^{-2}\, dx\, dy$) Maass form $u$ on $\bhv$
corresponds to $\frac1{\sqrt{2\p}} f_0$, where $f_0\in
L^2(\G\backslash G)$ has length~$1$ and weight~$0$. The factor $\p$
arises from the present normalization of the Haar measure on~$N$, and
the factor $2$ from $\G\backslash G \cong \left(\G\backslash
\bhv\right) \times \left( Z\backslash K\right)$, where $Z=\vz{I,-I}$.
\end{remark}

\begin{remark}\label{2-norm}
We note that density results entirely similar to those in this section
could have been obtained using the $2$-norm of the eigenvalues
$\ld_{\irr}$ in place of the $1$-norm.

Indeed, if one uses a zeta function involving $e^{-s \| {\ld_{\irr}
\|}_2}$ in place of $e^{-s \| {\ld_{\irr} \|}_1}$, the limit in
Theorem \ref{propv} is essentially the same, except for a different
multiplicative constant in the right-hand side. By using this limit
and following the arguments in this section, we obtain entirely
similar results with some changes in the multiplicative constants.
The proof of the analogue of Theorem \ref{propv} is very similar to
the one given in Sections 4-6, but one has to make a different choice
of the test functions, hence there are several computations and
estimates that need to be worked out again.
\end{remark}

\section{Sum Formula. Delta and Kloosterman terms}\label{prfpropv}
We shall use the sum formula in Theorem~2.7.1 of~\cite{BMP} in a way
similar to the application in Section~3 of that paper. For
 completeness, we will recall most of the notation. We refer the
 reader to~\cite{BMP} for any unexplained facts or notations.

 We apply the sum formula with $r=r'\in \oh'\setminus\vz0$. That
implies that the test functions~$k$ have the form $k=\times_{j=1}^d\,
k_j: \n \mapsto \prod_{j=1}^d k_j(\n_j)$ with all $k_j\in \Lkr$, see
Definition~\ref{Lkrdef}.

Throughout this section, we shall take the cusps~$\k$ and $\k'$ equal
to~$\infty$, and omit them from the notation.

The sum formula gives the following equality for each test
function~$k$:
\begin{equation}\label{sf*}
\int_Y k(\n) \, d\s_{r,r}(\n) = \Dt_{r,r}(k) + K_{-r,-r}(\B k).
\end{equation}
The integral on the left constitutes the \emph{spectral side} of the
sum formula. The measure $d\s_{r,r}$ contains information on the
spectral decomposition of $L^2(\G\backslash G)$ and on the Fourier
coefficients of the automorphic forms occurring in this
decomposition; see Sections~\ref{eisterm} and~\ref{discrterm}. The
\emph{geometric side} consists of the delta term $\Dt_{r,r}(k)$,
defined in~\eqref{Deltadef}, and the Kloosterman term $K_{-r,-r}(\B
k)$, see~\eqref{Kltdef}. The latter depends on a Bessel transform $\B
k$ of the test function~$k$.

In this section, we shall fix a special test function of product type,
leaving $k_j\in \Lkr$ free for $j \in E$ and choosing it in a special
way for $j \in Q$, depending on a parameter $s>0$. The purpose of
this section will be to investigate the behavior of the geometric
side as $s$ tends to~$0$. In \S\ref{Dtterm} we consider the delta
term. The study of the Kloosterman term, in \S\ref{Klterm} takes more
work. It turns out that the delta term gives the main contribution.

In all estimates, we take into account the dependence on $r$ and the
$k_j$, with $j\in E$.

The spectral side has the same behavior. This we shall use in
Section~\ref{discrterm} to prove Theorem~\ref{propv}.

\rmk{Notations}The map $\x\mapsto (\x^{\s_1},\ldots,\x^{\s_d})$ gives
an embedding of the number field~$F$ in~$\RR^d$.
We will often write $\x_j$ instead of $\x^{\s_j}$.

Accordingly, we define, for $x$, $y\in \RR^d$, the product
$xy\in\RR^d$ by $(xy)_j=x_jy_j$.

For $x\in\RR^d$, we put $S(x)\isdef \sum_{j=1}^d x_j$, extending the
trace $\Tr_{F/\QQ}$. Similarly, $N(y)\isdef \prod_{j=1}^d y_j$
extends the norm $N_{F/\QQ}$ to $N:(\RR^\ast)^d\rightarrow\RR^\ast$.

\rmk{Test functions}We write the set $\vz{1,\ldots,d}$ of places
of~$F$ as the disjoint union of three sets $E$, $Q_+$ and $Q_-$,
where $Q_+\cup Q_-\neq\emptyset$.

At the places $j\in E$, we keep $k_j\in \Lkr$ arbitrary. At the other
places, we make a special choice, depending on a parameter $s>0$:
\begin{alignat}2 \label{Gpchoice}
\text{If $j\in Q_+$:}&\quad &k_j(\n) &= \begin{cases}
e^{-s(1/4-\n^2)} &\text{ if }|\re\n|\leq \tau,\\
0&\text{ if }\n\in\frac12+\ZZ,\, |\n|>\tau;
\end{cases}\displaybreak[0]\\
\label{Gmchoice}
\text{if $j\in Q_-$:}&& k_j(\n)&= \begin{cases}0
&\text{ if }|\re\n|\leq \tau,\\
e^{-s(\n^2-1/4)}&\text{ if }\n\in\frac12+\ZZ,\, |\n|>\tau.
\end{cases}
\end{alignat}

\rmk{Norms of test functions} In Theorem~\ref{propv}, we have not
stated any uniformity of the limit in terms of the~$k_j$. In
Proposition~\ref{errorZeta}, we shall give some information on the
uniformity in~$k$. To do that, we now introduce some norms.

For $a>2$, let $\Lkr_a$ be the subspace of $k\in \Lkr$ for which
$k(\n) \ll \left(1+|\im\n|\right)^{-a}$ on the strip
$|\re\n|\leq\tau$.

For $\al\in [0,\tau]$ and $b\leq a$, we put:
\begin{align} \nonumber N_{\al,b}(k) &\isdef \sup_{\n:\, \re\n=\al}
\left(1+|\im\n|\right)^b |k(\n)| \text{ on } \Lkr_a,\\
\nonumber \Ndiscr(k) &\isdef \sum_{b\geq2,\, b\in2\ZZ}
{\textstyle\frac{b-1}2 \left| k\left(\frac{b-1}2\right)\right|}
\quad\text{on }\Lkr,
\\
\label{locnormdef*}
 \N_{\al,a}(k) &\isdef N_{0,a}(k) + N_{\al,a}(k) + \Ndiscr(k)
 \quad\text{on } \Lkr_a.
\end{align}
\noindent We extend $N_{\al,b}$ to~$\Lkr$ by defining it equal
to~$\infty$ outside~$\Lkr_a$.

We recall the definition of the integral transformation $\Eta(k)$
from~\eqref{Etadef}.
$$ \Eta (k) = \frac i2\int_{\re\n=0} k(\n)\, \n\,\tan\p\n\, d\n +
\sum_{b\geq2,\, b\in2\ZZ} {\textstyle\frac{b-1}2
k\left(\frac{b-1}2\right)}.
$$
Note that $\Eta$ is continuous on~$\Lkr_a$ with respect to
$N_{0,a}+\Ndiscr\leq \N_{\al ,a} $ for any $\al\in [0,\tau]$. We also
use the following notation:
\begin{equation}\label{normEdef} \|k\|_{\al,a,E} \isdef \prod_{j\in E}
\N_{\al,a}(k_j) \text{ if }k=\times_{j\in E}\, k_j \in \Lkr_a^{|E|}.
\end{equation}

\subsection{Delta term}\label{Dtterm}
Section~2.6 and Definition~2.5.2 in~\cite{BMP} give the definition of
the delta term:
\begin{equation}\label{Deltadef}
\Dt_{r,r}(k) = 2\vol(\G_N\backslash N) \prod_{j=1}^d\Eta(k_j).
\end{equation}

The group $N$ consists of the matrices $\matc 1x01$ with $x\in \RR^d$.
In~\cite{BMP}, we have chosen the Haar measure $dn =
\frac{dx_1}\p\cdots \frac{dx_d}\p$. The intersection $\G_N = \G\cap
N$ consists of the $\matc 1\x01$ with $\x\in \Okr\subset \RR^d$. So
the volume of $\G_N\backslash N$ is equal to $\p^{-d}\sqrt{|D_F|}$,
where $D_F$ is the discriminant of the number field~$F$; see, e.g.,
p.~115 in~\cite{La68}.

The factor $\al(r,r)$ in {\sl loc.\ cit.\ }is equal to $2$ in the
present context. To see this in Definition~2.6.1 in~\cite{BMP}, note
that the matrices $\matc\e00{1/\e}$ with $\e\in\Okr^\ast$ form a
system of representatives of $\G_N \backslash\G_P$. We take
$g_\infty=1$ at the cusp $\k=\infty$, hence $a_\g $ is equal to
$\matc{|\e^{\s_j}|}00{1/|\e^{\s_j}|}$ at place~$j$. So only $\e=\pm1$
contribute to~$\al(r,r)$, and $\ch_r(n_\infty(\g))=1$ for
$\g=\pm\matc1001$.

The factors $k_j$ of~$k$ with $j\in E$ are general. Let us put
$\Eta_E(k) \isdef \prod_{j\in E} \Eta(k_j)$.

We are left with the factors for $j\in Q$. We have chosen the
corresponding~$k_j$ in~\eqref{Gpchoice} and~\eqref{Gmchoice}. For
$j\in Q_+$:
\begin{align*}
\Eta(k_j) &= \frac i2 \int_{\re\n=0} e^{s(\n^2-1/4)} \n\,\tan\p\n\,
d\n + \tfrac12 \\
&= s^{-1} e^{-s/4} \int_{t=0}^\infty e^{-t^2} t \, dt + \oh\left(
\int_0^\infty e^{-st^2-s/4} t e^{-2\p t}\, dt \right)+ \tfrac12\\
&= \tfrac12 s^{-1} + \oh(1),
\end{align*}
and for $j\in Q_-$:
\begin{align*}
\Eta(k_j) &= \sum_{m=2}^\infty\left( m-\tfrac12\right) e^{s(m-m^2)}.
\end{align*}
In order to estimate this quantity, we replace the sum by
$\int_1^\infty f(x)\,dx$, with $f(x)=(x-\nobreak\frac12) \allowbreak
e^{-s (x^2-x)}$, and then we need to estimate the error.

 The maximum of $f$ occurs at $x_{\rm max}=\frac12+\frac1{\sqrt{2s}}$.
We have $\int_1^\infty f(x)\,dx=\tfrac12 s^{-1}+ \oh(1)$. To estimate
the error we consider two intervals: The sum over $2\leq m\leq
m_1=\left[ x_{\rm max}\right]$ is larger than the integral over
$[1,m_1]$, and the difference between sum and integral is smaller
than $f(x_{\rm max})$. The sum over $m\geq m_1+1$ is larger than the
integral over $[m_1+1,\infty)$, and again the difference is smaller
than $f(x_{\rm max})$. The missing integral over $[m_1,m_1+\nobreak
1]$ is smaller than $f(x_{\rm max})$ as well. So the error is $\ll
f(x_{\rm max}) \ll \frac1{\sqrt{s }} $ and we obtain
\begin{align}
\nonumber \Eta(k_j) &= \tfrac12 s^{-1} + \oh\left( s ^{-1/2} \right),
\\
\label{Dtest}
\Dt_{r,r}(k) &= \frac{2^{1+|E|}}{(2\p)^d} \sqrt{|D_F|}\, \Eta_E(k)\,
s^{|E|-d}\\
\nonumber
&\qquad\qquad\hbox{} \cdot \left( 1+\,(\text{if $Q_-\neq\emptyset$})\,
\oh (s^{1/2}) +\oh (s) \right).
\end{align}
With the convention in~\eqref{normEdef}, we can restate this as
\begin{prop}\label{Dtest1} Let $k=\times_{j=1}^d k_j$, with $k_j\in
\Lkr$ arbitrary for $j\in E$, and with $k_j$ as defined in
\eqref{Gpchoice} and \eqref{Gmchoice}, for $j\in Q$. Then we have for
each $\al\in [0,\tau]$, as $s\downarrow0$:
\begin{equation*}
 \Dt_{r,r}(k) = \frac{2^{1+|E|}}{(2\p)^d} \sqrt{|D_F|}\,
 \Eta_E(k)\,s^{|E|-d} + \|k\|_{\al,a,E}\,\oh_F(s^{|E|-d + c}).
\end{equation*}
where $c=\frac 12$ if $Q_-\ne \emptyset$ and $c=1$ if $Q_-=\emptyset$.
\end{prop}

\subsection{Kloosterman term}\label{Klterm}
In \S3.3 of~\cite{BMP} it was sufficient for our purposes to estimate
the Kloosterman term by a quantity that is of the same order as the
delta term. Here we want to use the main term in~\eqref{Dtest} for an
asymptotic result. So we have to do better on the Kloosterman term.

We see no way of using interference between Kloosterman sums, so we
shall estimate all terms by their absolute value, employing a
Sali\'e-Weil type estimate of individual Kloosterman sums.

For any function $f:(\RR^\ast)^d\rightarrow \CC$ that decreases
sufficiently fast as the $y_j$ tend to zero and infinity, we define
the following sum of Kloosterman sums:
\begin{equation}\label{Kltdef}
K_{r,r}(f) \isdef \sum_{c\in \qq,\, c\neq0} \frac{S(r,r;c)}{|N(c)|}
f\left( \frac{r^2}{c^2} \right).
\end{equation}
Here $\qq$ is the ideal $\qq\subset\Okr$ such that $\G=\G_0(\qq)$. We
have $\left( \frac {r^2}{c^2} \right)_j = \left( \frac {r^2}{c^2}
\right)^{\s_j}$. The quantity $|N(c)|$ is the norm of the ideal
$(c)\subset\Okr$. The Kloosterman sum is defined by
\begin{equation} \label{Klsumdef} S(r,r;c) =
\mathop{\textstyle{\sum^\ast}}_{d \bmod\, c} e^{2\p
i\Tr_{F/\QQ}(r(d+a)/c)},
\end{equation}
where $d$ runs over representatives of $\Okr\bmod\, (c)$ for which
there exists $a\in\Okr$ such that $ad\equiv 1\, \bmod\, (c)$.

The aim of this subsection is to prove the following estimate for the
Kloosterman term.
\begin{prop}\label{KLTest} Let $k=\times_{j=1}^d k_j$, with $k_j\in
\Lkr$ arbitrary for $j\in E$, and with $k_j$ as defined in
\eqref{Gpchoice} and \eqref{Gmchoice} for $j\in Q$. For $\al\in
\left(\frac12,\tau\right]$, and $0<\e<1-\tau$, as $s\downarrow0$, we
have:
\begin{equation*}
 K_{-r,-r}(\B k) \le C(\qq,r,\e)\, \|k\|_{\al,a,E} \,
 s^{-(3/4+\e)|Q_+|-(1/4+1/8\al+\e)|Q_-| },
\end{equation*}
with $C(\qq,r,\e)=N(\qq)^{-1/2+\e} |N(r)|^{1/2+\e}$.
\end{prop}

The choice $\al=\tau$ is optimal in the $s$-aspect. We see that the
Kloosterman term is indeed of smaller order than
$s^{|E|-d}=s^{-|Q_+|-|Q_-|}$ in~\eqref{Dtest}.

\rmk{Bessel transform} The transform $\B k$ occurring in the
Kloosterman term is of product type: $\B k(y) = \prod_{j=1}^d \bt_+
k_j(y_j)$ for $y\in (0,\infty)^d$, with the following Bessel
transformation:
\begin{align}\label{btplusdef}
 \bt_+ k(y) &= \frac i2\int_{\re\n=0} k(\n) \left( J_{-2\n}(4\p\sqrt
y) - J_{2\n}(4\p\sqrt y)\right) \frac{\n\,d\n}{\cos\p\n}\\
\nonumber &\qquad\hbox{} + 2 \sum_{b\geq 2,\, b\in 2\ZZ} (-1)^{b/2}
k{\textstyle\left(\frac{b-1}2\right) \frac{b-1}2} J_{b-1}(4\p\sqrt
y), \displaybreak[0] \\
\nonumber J_w(t) &= \sum_{n=0}^\infty \frac{(-1)^n}{n!\, \Gf(w+n+1)}
\left(\frac t2\right)^{w+2n}.
\end{align}
In (25) and~(26) of~\cite{BMP}, we have rewritten the integral
defining $\bt_+k$ in several ways. As in \S3.3 of~\cite{BMP}, we use
some estimates for the Bessel function to find bounds for $\bt_+k$.
With $0\leq \al \leq \tau$:
\begin{alignat}2 \label{Jestsmally}
&\text{For }\re\n=\al&\quad J_{2\n}(y) &\ll y^{2\al}e^{\p|\im\n|}
\left(1+|\im\n|\right)^{-2\al-1/2}\\
\nonumber
&&&\qquad\text{ as }y\downarrow0, \pagebreak[0] \\
\nonumber
&\text{for }b\in2\ZZ,\, b\geq 2&\quad J_{b-1}(y) &\ll
y^{b-1}\Gf(b)^{-1}\text{ as }y\downarrow0, \pagebreak[0] \\
\nonumber
&\text{for }\re\n=0& J_{2\n}(y) &\ll e^{\p|\im\n|}\text{ for all }y>0,
\pagebreak[0] \\
\nonumber
&\text{for }u>0& J_u(y) &\ll u^{-1}\text{ for all }y>0.
\end{alignat}
For convenience, we write $t=4\p\sqrt y$. We fix $\al
\in\left(\frac12,\tau\right]$. For $k\in\Lkr_a$, with $a>2$, we find
the following estimates along the same lines as in~\cite{BMP}, \S3.3.
\begin{align*}
\bt_+k(y) &= -i\int_{\re\n=\al} k(\n) J_{2\n}(t)
\frac{\n\,d\n}{\cos\p\n} + 2\sum_{b\geq4} (-1)^{b/2}{\textstyle
\frac{b-1}2 k\left(\frac{b-1}2\right)} J_{b-1}(t)\\
&\ll_\al N_{\al,a}(k) y^\al + y^{3/2} \Ndiscr(k) \ll \N_{\al,a}(k)
y^\al
\quad\text{ as }y\downarrow0, \pagebreak[0] \\
\bt_+k(y) &= -i\int_{\re\n=0} k(\n) J_{2\n}(t)
\frac{\n\,d\n}{\cos\p\n} + 2 \sum_{b\geq2} (-1)^{b/2} {\textstyle
\frac{b-1}2 k\left(\frac{b-1}2\right)} J_{b-1}(t)\\
&\ll N_{0,a}(k) + \Ndiscr(k) \ll \N_{\al,a}(k)\quad\text{ as }
y\rightarrow\infty, \pagebreak[0]
\\
\bt_+k(y) &\ll \N_{\al,a}(k) \min\left(y^\al,1\right)\quad\text{ for
}y>0.
\end{align*}

If we specialize $k=k_j$ as in~\eqref{Gpchoice} and~\eqref{Gmchoice}
for $j\in Q$, we can obtain better estimates by reconsidering the
integrals. First we note that for $\re\n>\alpha$:
\begin{equation}
J_{2\n}(y) = \frac1{2\p i} \int_{\re w=-\al} \left( \frac
y2\right)^{-2w} \frac{\G(\n+w)}{\G(1+\n-w)}\, dw.
\end{equation}
(To derive this integral representation from the power series
expansion, move the line of integration to the left.) We take
$\frac12<\al<\g<\al+\frac12$, and write $\n=\g+iq$, $w=-\al+it$. To
estimate $J_{2\n}$, it suffices to consider the case $q\geq0$. We
find for all $y>0$:
\begin{align*}
J_{2\n}&(y) \ll_{\al,\g} y^{2\al} \int_{-\infty}^\infty e^{-\frac\p2|
q+t|+\frac\p2| q-t|} (1+| q+t|)^{\g-\al-1/2} \\
&\qquad\qquad \hbox{} \cdot (1+| q-t|)^{-\g-\al-1/2}\, dt
\displaybreak[0] \\
&= y^{2\al} e^{\p q} \int_0^\infty (1+x)^{\g-\al-1/2} (1+2
q+x)^{-\g-\al-1/2}\, dx\\
&\qquad\hbox{} + y^{2\al} q \int_{-1}^1 e^{-\p q x} (1+ q(1+
x))^{\g-\al-1/2}(1+ q(1- x))^{-\g-\al-1/2}\, d x\\
&\qquad\hbox{} + y^{2\al} e^{-\p q} \int_0^\infty (1+x)^{-\g-\al-1/2}
(1+2 q+x)^{\g-\al-1/2}\, dx \displaybreak[0] \\
&\ll y^{2\al} e^{\p q}\left( \int_1^{ q+1} x^{\g-\al-1/2} (1+
q)^{-\g-\al-1/2}\, d x + \int_{ q+1}^\infty x^{-2\al-1}\, d x\right)
\\
&\qquad\hbox{} + y^{2\al} q \left( \int_0^1 e^{\p q x} (1+
q)^{-\al-\g-1/2}\, d x + (1+ q)^{\g-\al-1/2}\right)
\\
&\qquad\hbox{} + y^{2\al} e^{-\p q} q^\g\int_0^\infty
(1+x)^{-2\al-1}\, dx \displaybreak[0] \\
&\ll y^{2\al}e^{\p q} (1+ q)^{1/2-\al-\g}.
\end{align*}
 The best choice of $\g$ seems $\al+\frac12-\e$, with $\e>0$ small.
We'll apply the estimate with $\frac12<\g<\frac32$, hence we take
$\e<\frac12$. The advantage of the present estimate
above~\eqref{Jestsmally} is its validity for all $y>0$.

In the case of real~$\n$, $\n\geq \frac32$, we find
\begin{align*}
J_{2\n}&(y) \ll_\al y^{2\al} \int_{-\infty}^\infty
e^{-t(\arg(\n-\al+it)-\arg(1+\n-\al+it))}
\frac{|\n-\al+it|^{\n-\al-1/2}}{|1+\n+\al|^{\n+\al+1/2}}\,dt
\displaybreak[0]
\\
&\ll y^{2\al} \int_{-\infty}^\infty |1+\n+\al+it|^{-2\al-1}\, dt \ll
y^{2\al} \n^{-2\al}.
\end{align*}\smallskip

We apply these estimates to find bounds for the local Bessel
transforms. We take $\al\in\left(\frac12,\tau\right]$,
$\g=\al+\frac12-\e$, $0<\e<\frac12$. For $j\in Q$, the functions
$k_j$ are holomorphic on $|\re\n| < \frac32$, so $k_j(\n)$ makes
sense if $\re\n=\g$.

For $j\in Q_+$, we have uniformly for $0<s<1$:
\begin{align*}
\bt_+&k_j(y) = -i \int_{\re\n=\g} k_j(\n) J_{2\n}(t) \frac{\n\,
d\n}{\cos\p\n} \displaybreak[0] \\
&\ll \int_{-\infty}^\infty e^{s(\g^2-1/4-u^2)} y^\al
(1+|u|)^{1/2-\al-\g}|u|\, du \displaybreak[0] \\
&\ll y^\al s^{\al/2+\g/2-5/4} \int_0^\infty e^{-u^2} (\sqrt
s+u)^{1/2-\al-\g}u\, du \displaybreak[0] \\
&\ll y^\al s^{\al-1-\e}, \displaybreak[0] \\
\bt_+ k_j(y) &\ll e^{-s/4} \int_{-\infty}^\infty e^{-sw^2}
\frac{|w|e^{\p w}}{\cosh\p w}\, dw+1\displaybreak[0] \\
&\ll e^{-s /4} \int_0^\infty e^{-s w^2} w\, dw + 1\\
&\ll s^{-1}, \displaybreak[0] \\
\bt_+k_j(y) &\ll \min\left( s ^{\al-1-\e} y^\al, s^{-1} \right)
\quad\text{ for }y>0.
\end{align*}
In the last estimate, we have the parameter $s$ inside the minimum.
This will enable us to improve the estimates in~\cite{BMP}, \S3.3.

For $j\in Q_-$ we find a similar, but slightly better estimate,
uniformly for $0<s<1$:
\begin{align*}
\bt_+k_j(y) &= 2 \sum_{b\geq 4,\, b\in 2\ZZ} (-1)^{b/2} \tfrac{(b-1)}2
e^{s \frac b2(1-\frac b2)} J_{b-1}(t) \displaybreak[0] \\
&\ll \sum_{b\geq 4,\, b\in 2\ZZ} y^\al (b-1)^{1-2\al}
e^{-s((b-1)/2)^2} \displaybreak[0] \\
&\ll y^\al \left( \left(\sqrt{\frac{2-2\al}{2s}}\right)^{1+1-2\al} 1 +
\int_{x=\sqrt{\frac{2-2\al}{2s}}}^\infty x^{1-2\al} e^{-sx^2}\, dx
\right) \displaybreak[0] \\
&\ll_\al y^\al s^{\al-1}, \displaybreak[0] \\
\bt_+ k_j(y) &\ll \sum_{b\geq 4,\, b\in 2\ZZ} e^{s \frac b2(1-\frac
b2)} \leq \int_{x=1}^\infty e^{s x(1-x)} \, dx\ll s^{-1/2},
\displaybreak[0] \\
\bt_+ k_j(y) &\ll \min\left( y^\al s^{\al-1}, s^{-1/2} \right)
\qquad\text{ for all }y>0.
\end{align*}

\noindent In this way we have proved:
\begin{lem}\label{Besselest} Let $\frac12<\al\leq \tau$,
$0<\e<\frac12$. For $k=\times_{j=1}^d k_j$, with $k_j\in \Lkr$
arbitrary for the places $j\in E$, and $k_j$ as in \eqref{Gpchoice}
and \eqref{Gmchoice} for $j\in Q$, we have uniformly for $s\in
(0,1)$:
\begin{align*}
 \bt_+k_j(y) &\ll \N_{\al,a}(k) \min\left(y^\al,1\right),\quad
 \text{for $j\in E$,}\\
  \bt_+k_j(y) &\ll \min\left(s^{\al-1-\e} y^{\al} , s^{-1} \right),
\quad \text{for $j\in Q_+$}\\
  \bt_+k_j(y) &\ll \min\left( s^{\al-1} y^{\al} , s^{-1/2} \right),
\quad \text{for $j\in Q_-$}.
\end{align*}
\end{lem}

\rmk{Sum over the units}We shall apply Lemma~8.1 in~\cite{BMste} in
the following form:
\begin{lem}\label{su}Let $a,b\in\RR$, $a+b>0$. Let $p_j, q_j>0$ for
$j=1,\ldots,d$. There exists $C\geq0$ such that for all
$f:(\RR^\ast)^d\rightarrow\CC$ satisfying
\[ |f(y)|\leq \prod_{j=1}^d\min\left(p_j|y_j|^a,q_j|y_j|^{-b}\right),
\]
we have
\begin{align*} \sum_{\e\in\Okr^\ast} \left|f(\e y) \right| &\ll
\min\left( N(p)\,|N(y)|^a, N(q)\, |N(y)|^{-b} \right)\\
&\qquad\hbox{} \cdot \left( 1 + \left| \log|N(y)| +\frac1{a+b}
\log\frac{N(p)}{N(q)} \right|^{d-1} \right).
\end{align*}
\end{lem}
\begin{proof}This is a direct consequence of Lemma~8.1
in~\cite{BMste}. There we had to take into account the complex places
of~$F$. We apply the lemma with $e=0$ and all $n_j=1$.

With $\eta_j = p_j^{1/(a+b)}q_j^{-1/(a+b)}|y_j|$, we have
\[ \min\left(p_j|y_j|^a,q_j|y_j|^{-b}\right) = p_j^{b/(a+b)}
q_j^{a/(a+b)} \min\left(\eta_j^a,\eta_j^{-b}\right).\]
Lemma~8.1 in~\cite{BMste} bounds the sum over the units by
\begin{align*} &N(p)^{b/(a+b)} N(q)^{a/(a+b)} \left( 1+\left|\log
N(\eta)\right|^{d-1} \right) \min\left( N(\eta)^a,N(\eta)^{-b}
\right)\\
&\ll \left( 1 + \left|\log\left|\vphantom{\frac12}\right.
\frac{N(p)^{1/(a+b)}}{N(q)^{1/(a+b)}}N(y)
\left.\vphantom{\frac12}\right| \right|^{d-1} \right) \min\left( N(p)
|N(y)|^a, N(q) |N(y)|^{-b} \right).
\end{align*}
\end{proof}

\rmk{Kloosterman term} We will use the fact that the Kloosterman sums
given in \eqref{Klsumdef} satisfy an estimate of Weil-Sali\'e type
(see Theorem~10 in~\cite{BM95}):
\begin{equation}\label{SWE}
S(r,r;c) \ll_{F,\e} N_{r,r}(c)^{1/2} |N(c)|^{1/2+\e},
\end{equation}
for each $\e>0$. The factor $N_{r,r}(c)$ is described in terms of the
following decompositions in prime ideals:
\[ (r) = \prod_P P^{v_P(r)},\qquad (c) = \prod_P P^{v_P(c)},\qquad
\Okr' = \prod_P P^{-d_P}.\]
Then $N_{r,r}(c) = \prod_P N(P)^{\min(v_P(r),v_P(c)-d_P)}$.

For all test functions~$k$ in the sum formula, the sum $K_{r,r}(f)$
converges absolutely. This convergence is part of the statement of
Theorem~2.7.1 in~\cite{BMP}. For our choice of test functions, the
absolute convergence follows from the estimates of Bessel transforms
obtained above, and the next lemma. It is a slight improvement of
Lemma~3.2.1 in~\cite{BMP}.

\begin{lem}\label{ks} Let $f:(\RR^\ast)^d\rightarrow\CC$ satisfy
\[ |f(y) | \leq \prod_{j=1}^d \min\left(p_j(f),q_j(f)|y_j|^\al\right),
\]
with $p_j(f),\, q_j(f)>0$, for $j=1,\ldots,d$, and with $\al>1/4$. Let
$r\in\Okr'\setminus\vz0$.

\noindent Then $K_{r,r}(f)$ converges absolutely, and for each $\e>0$,
 sufficiently small:
\[ K_{r,r}(f) \ll_{F,\e} N(p_f) \left(
\frac{N(q_f)}{N(p_f)}\right)^{1/4\al+\e} N(\qq)^{-1/2+\e}
|N(r)|^{1/2+\e}, \]
with $N(p_f)=\prod_{j=1}^d p_j(f)$, and $N(q_f)$ similarly.
\end{lem}
\begin{proof}We follow closely the proof of Lemma~3.2.1 in~\cite{BMP}.
\begin{align*}
|K_{r,r}(f)| &\leq \sum_c \frac{|S(r,r;c)|}{|N(c)|} \prod_{j=1}^d
\min\left( p_j(f), q_j(f) |r_j|^{2\al} |c_j|^{-2\al}\right)\\
&\ll_{F,\e} \sum_{(c)\subset \qq,\, (c)\neq(0)} N_{r,r}(c)^{1/2}
|N(c)|^{\e-1/2}\\
&\qquad\hbox{} \cdot \sum_{\z\in\Okr^\ast} \prod_{j=1}^d \min\left(
 p_j(f), q_j(f) |r_j|^{2\al} |\z_jc_j|^{-2\al}\right).
\end{align*}
We use Lemma~\ref{su} to estimate the sum over~$\z\in\Okr^\ast$, where
we take $a=0$, $b=2\al$, $p_j=p_j(f)$, $q_j=q_j(f)$, and $y=c/r$.
This gives the bound
\begin{align*}
\sum_\z \ast &\ll_{F,\e,\al} \min\left( N(p_f), N(q_f)
\left|N(r/c)\right|^{2\al} \right)
\\
&\hbox{} \cdot \left( 1 + \left| \log|N(c/r)| +
\frac1{2\al}\log(N(p_f)/N(q_f)) \right|^{d-1} \right).
\end{align*}

\noindent As in the proof of Lemma~3.2.1 in~\cite{BMP}, we write $(r)
= R_+R_-^{-1}$, with $R_+=\prod_{P,\,v_P(r)\geq0}P^{v_P(r)}$,
$R_-=\prod_{P,\, v_P(r)<0}P^{-v_P(r)}$; hence $N(R_-)\leq
N(\Okr')^{-1} \ll_F1$. Any ideal $(c)$ in the remaining sum is
written uniquely as $(c)=KJ$, where $K=((c),R_+)$, and hence
$L=R_+K^{-1}$ and $J=(c)K^{-1}$ are prime to each other. Let us write
$\qq=\qq_1\qq_2$, with $\qq_1=(\qq,R_+)$ and $\qq_2$ relatively prime
to $R_+\qq_1^{-1}$. As we sum over $(c)\subset\qq$, we obtain the
additional conditions $\qq_1\supset K$ and $\qq_2\supset J$. We
proceed as in the proof of Lemma~3.2.1 of~\cite{BMP}, and obtain:
\begin{align*}\left|K_{r,r}(f)\right|&\ll \sum_{K,\, \qq_1|K| R_+}
\sum_{J,\, (J,L)=1,\, \qq_2|J} \frac{ N(K)^\e}{N(J)^{1/2-\e}}
\min\left( N(p_f),N(q_f) \frac{N(L)^{2\al}}{N(J)^{2\al}}\right)
\\
&\qquad\hbox{} \cdot \left( 1 + \left| \log \frac
{N(J)N(p_f)^{1/2\al}}{N(L)N(q_f)^{1/2\al} } + \log N(R_-)
\right|^{d-1} \right) \displaybreak[0]\\
&= \sum_{K,\, \qq_1|K|R_+} S_K.
\end{align*}
 The inner sum $S_K$ can be majorized by the same sum taken over all
ideals $J$ such that $\qq_2 | J$. We use that the number the ideals
in $\Okr$ with norm $n$ is $\oh(n^\e)$ and obtain:
\begin{align*}
S_K &\ll N(K)^\e \sum_{n=1}^\infty \left(n N(\qq_2) \right)^{2\e-1/2}
N(p_f) \min\left(1, \left(\frac{\k N(L)}{n }\right)^{2\al} \right)\\
&\quad\hbox{} \cdot \left( \oh_F(1) + \left| \log \frac{n
}{N(L)\k}\right|^{d-1} \right) ,
\end{align*}
with $\k=N(q_f)^{1/2\al} N(p_f)^{-1/2\al}N(\qq_2)^{-1}$. We split up
the sum at $n \approx N(L)\k$ and we use that $L=R_+K^{-1}$ .

The sum over $n\leq N(L)\k$ is present only if $N(L)\geq \k^{-1}$.
Under this condition, or equivalently, for $N(K)\leq \k N(R_+)$, we
obtain:
\begin{align}
\label{kaplarge}
S_K&\ll N(p_f) N(\qq_2)^{-1/2+2\e} N(K)^\e \left(
N(L)\k\right)^{1/2+3\e}
\\
\nonumber
&\quad\hbox{} + N(p_f) N(\qq_2)^{-1/2+2\e} N(K)^\e
(N(L)\k)^{1/2+3\e-2\al} \k^{2\al} N(L)^{2\al} \pagebreak[0] \\
\nonumber
&\ll N(p_f) N(\qq_2)^{-1/2+2\e} N(K)^{-1/2+2\e} N(R_+)^{1/2+3\e}
\k^{1/2+3\e} . \end{align}

If $N(L)<\k^{-1}$, or equivalently, for $N(K)> \k N(R_+)$, we find
\begin{align}
\label{kapsmall}
S_K &\ll N(p_f) N(\qq_2)^{-1/2+2\e} N(K)^\e \left(\k
N(L)\right)^{2\al-\e}
\\
\nonumber
&\ll N(p_f) N(\qq_2)^{-1/2+2\e} N(K)^{-2\al+2\e} N(R_+)^{2\al-\e}
\k^{2\al-\e} . \end{align}

We look for the maximum of these estimates for $K|R_+$, $\qq_1|K$. If
$\k\geq 1$ we need solely \eqref{kaplarge}, and find
\[ S_K \ll N(p_f) N(\qq_1)^{-1/2+2\e} N(\qq_2)^{-1/2+2\e}
N(R_+)^{1/2+3\e} \k^{1/2+3\e} . \]
If $N(\qq_1)N(R_+)^{-1}\leq \k <1$:
\begin{align*}
S_K &\ll \max\left( N(p_f) N(\qq_1)^{-1/2+2\e} N(\qq_2)^{-1/2+2\e}
N(R_+)^{1/2+3\e} \k^{1/2+3\e},\right.
\\
&\qquad \qquad\qquad\left. N(p_f) N(\qq_2)^{-1/2+2\e} N(R_+)^\e \k^\e
\right) \pagebreak[0]
\\
&\ll N(p_f) N(\qq_1)^{-1/2+2\e} N(\qq_2)^{-1/2+2\e} N(R_+)^{1/2+3\e}
\k^{1/2+3\e}.
\end{align*}
Finally, if $\k<N(\qq_1)N(R_+)^{-1}$:
\begin{align*}
 S_K &\ll N(p_f) N(\qq_1)^{-2\al+2\e} N(\qq_2)^{-1/2+2\e}
 N(R_+)^{2\al-\e} \k^{2\al-\e} \\
&\ll N(p_f) N(\qq_1)^{-1/2+2\e} N(\qq_2)^{-1/2+2\e} N(R_+)^{1/2+3\e}
\k^{1/2+3\e}.
\end{align*}
The number of possibilities of choosing $K$ is $\oh\left( N(R_+)^\e
N(\qq_1)^{-\e}\right)$. We find, for any $\k$,
\begin{align*}
K_{r,r}(f) &\ll N(p_f) N(\qq_2)^{-1/2+2\e} N(\qq_1)^{-1/2+2\e}
N(R_+)^{1/2+4\e} \k^{1/2+3\e} . \end{align*}

\noindent
We use $N(R_+) \ll_F |N(r)|$ to complete the proof.

We observe that $\al>\frac14$ is indeed necessary for this proof to go
through. We need $2\al+\frac12-\e>1$ to have convergence in the
 series leading to \eqref{kaplarge} and~\eqref{kapsmall}, and we need
$2\al-\frac12-4\e>0$ in a later estimate (case $\k<
N(\qq_1)N(R_+)^{-1}$).
\end{proof}
\smallskip

Now we are in a position to prove Proposition \ref{KLTest}. Let $f=\B
k$ and let $p_j(f) =s ^{-1}$ (resp. $p_j(f) =s ^{-1/2}$) if $j \in
Q_+$ (resp. $Q_-$). Furthermore, let $q_j(f) = s ^{\al-1-\e}$ (resp.
$q_j(f) = s ^{\al-1}$), if $j\in Q_+$ (resp. $Q_-$). Let
$p_j(f)=q_j(f)=\N_{\al,a}(k_j)$, if $j\in E$. By Lemma
\ref{Besselest} we have
\[ |f(y) | \leq \prod_{j=1}^d \min\left(p_j(f),q_j(f)|y_j|^\al\right).
\]
We also note that
\begin{align*}N(p_f)&=\|k\|_{\al,a,E} s^{-(|Q_+|+\frac 12|Q_-|)} ,\\
 N(q_f)&=\|k\|_{\al,a,E} s^{-((1+\e - \al)|Q_+| +(1-\al)|Q_-|)}.
\end{align*}
Now, if $C(\qq,r,\e) = N(\qq)^{-1/2+\e} |N(r)|^{1/2+\e}$, Lemma
 \ref{ks} and a calculation show that
\begin{equation}
K_{-r,-r}(\B k) \ll \|k\|_{\al,a,E} C(\qq,r,\e) s^{-((3/4+\e)|Q_+| +
(1/4+1/8\al+\e)|Q_-|)}
\end{equation}

\section{Sum Formula. The Eisenstein term}\label{eisterm}

The goal of this section is to give an estimate for the Eisenstein
term of the sum formula. We start by introducing some notation.

The set $Y$ in the sum formula~\eqref{sf*} has the form
\[Y \isdef \left( i[0,\infty) \cup (0,{\textstyle\frac12}) \cup \vzm{
{\textstyle\frac{b-1}2}} {b\geq2, \, b\in 2\ZZ} \vphantom{R^{R^4}}
\right)^d.\]
The measure $d\s_{r,r}(\n)$ on~$Y$ in the left hand side
of~\eqref{sf*} is the sum of a measure $d\s_{r,r}^{\rm discr}$, to be
 discussed in Section~\ref{discrterm}, and a measure $d\s_{r,r}^{\rm
 cont}$, given by
\begin{align} \label{measure}
\int_Y f(\n) \, &d\s^{\rm cont}_{r,r}(\n) \isdef \sum_{\k\in \Pkr}
c_\k \sum_{\m \in {\mathcal L}_\lambda} \int_{-\infty}^\infty
f(iy+i\m)\, \left| D_\k^{\infty,r}(iy,i\m) \right|^2 \, dy
\end{align}
for compactly supported continuous functions on~$Y$. All test
functions $k=\times_{j=1}^d k_j$, $k_j\in \Lkr$, are integrable for
this measure, in particular if we take $k_j$ as in \eqref{Gpchoice}
and~\eqref{Gmchoice}.

$\Pkr$ is a full set of representatives of $\G_0(\qq)$-inequivalent
cusps. For each $\k \in {\mathcal P}$, the number $c_\k$ is positive;
$D_\k^{\infty,r}(\n,i\m)$ is a normalized Fourier coefficient of the
Eisenstein series $E_q(P^\k,\n,i\m,g)$ with weight $q$ and spectral
parameter $\n\in \CC$, and $\mu$ is in a lattice $\Lkr_\k$ in the
hyperplane $\sum_{j=1}^d x_j=0$ depending on~$\k$. \smallskip

 To show that the contribution of the Eisenstein term is negligible in
the context of this paper, we need an estimate for
$D_\k^{\infty,r}(\n,i\m)$ in terms of $\im \n$ and $\mu$. This will
take quite some work, which we carry out in this section. It is at
this point that we essentially use that we are in an arithmetic
situation. In general, the contribution of the Eisenstein series
could be comparable to, or even larger than that of the discrete
spectrum. Here, there are expressions for the Fourier coefficients of
Eisenstein series in terms of number theoretic quantities, for which
one has additional information.

\subsection{Fourier coefficients of Eisenstein series}\label{FcE} In
view of (18) in \cite{BMP}, in order to estimate the coefficient
$D_\k^{\infty,r}(\n,i\m)$, it suffices to consider the weight $q=0\in
(2\ZZ)^d$. As usual, let $\G(\qq) \subset \G_0 (\qq)$ be the
principal congruence subgroup of level $\qq$, where $\qq$ is an ideal
in $\Okr$.

We denote $a[y]=\vz{ \matc{\sqrt {y_1}}00{1/\sqrt{y_1}},\cdots,
\matc{\sqrt {y_d}}00{1/\sqrt{y_d}}} \in G$ for $y\in \RR_{>0}^d$, and
$a[y]^\rho=\prod_j y_j^{1/2}$, $a[y]^{i\m }=\prod_j y_j^{i\m_j}$. The
Eisenstein series is, for $\re\n>\frac12$, given by the sum
\begin{equation}
E_0(P^\k,\n,i\m;g)=\sum_{\g\in \G_0(\qq)_\k \backslash \G_0(\qq)}
a_\k(\g g)^{\r+2\n\r+i\m},
\end{equation}
where $g_\k \in G$ satisfies $\k=g_\k \infty$, where $g=g_\k n_\k(g)
a_\k(g) k_\k(g)$ for each $g\in G$, with $n_\k(g) \in N=\vz{\matc
1\ast01}$, $a_\k(g) \in \vz{a[y]}{y\in \RR_{>0}^d}$, $k_\k(g)\in K$,
and where $\G_0(\qq)_\k$ is the subgroup fixing the cusp~$\k$. The
lattice $\Lkr_\k$ consists of the $\m$ in the hyperplane
$\sum_{j=1}^dx_j=0$ that satisfy $a_\k(\dt)^{i\m}=1$ for all $\dt\in
\G_0(\qq)_\k$.

The Eisenstein series $E_0(P^\k,\n,i\m,g)$ is a linear combination of
Eisenstein series $E_0^{\G(\qq)}(P^{\k_1},\n,i\m,g)$ for the
principal congruence subgroup $\G(\qq)$, with $\k_1$ running through
the cusps of $\Gamma(\qq)$ above $\k$. The coefficients in this
linear combination depend on the choice of the $g_{\k_1} \in \GG_\QQ$
such that $\k_1 =g_{\k_1} \infty$. This choice causes factors of the
form $a[p]^{\r+2\n \r + i\m}$, so it will not influence the size if
$\re\, \nu =0$.

To see this, we note that if $\g$ runs through a set of
representatives of $\G_0(\qq)_\k\backslash\G_0(\qq)$, then $\g^{-1}$
runs through a set $R\subset \G_0(\qq)$ such that each element $\g\k$
of the orbit $\G_0(\qq)\k$ occurs exactly once. This orbit consists
of a finite number of $\G(\qq)$-orbits, for which we choose
representatives $\g_j\in \G_0(\qq)$. Let $\k_j=\g_j \k$. So there are
subsets $R_j\subset \G(\qq)$ such that $R = \bigsqcup_j R_j\g_j$.

The $\dt\in R_j$ represent the orbit $\G(\qq)\k_j$. It is not hard to
check that this implies that $\dt^{-1}$ runs through representatives
of $\G(\qq)_{\k_j}\backslash \G(\qq)$ when $\dt$ runs through $R_j$.

Let us choose $g_{\k_j} = \g_j g_\k$. Then
\begin{align*}
g_\k &n_\k( \g_j^{-1}\dt^{-1}g) a_\k( \g_j^{-1}\dt^{-1}g) k_\k(
\g_j^{-1}\dt^{-1}g) = \g_j^{-1}\dt^{-1} g\\
& = \g_j^{-1} g_{\k_j} n_{\k_1}( \dt^{-1}g) a_{\k_1}(\dt^{-1}g)
k_{\k_1}(\dt^{-1}g)
\end{align*}
shows that $a_{\k_1}(\dt^{-1}g) = a_\k(\g_j^{-1}\dt^{-1}g)$, and hence
\[ E_0(P^\k,\n,i\m;g) = \sum_j E^{\G(\qq)}_0(P^{\k_j},\n,i\m;g).\]
Any other choice of $g_{\k_1}$ has the form $g_{\k_1}=\g_j g_\k pl$,
with $p\in NA$, $l=\pm I$. This will indeed cause a factor
$a[p]^{\r+2\n\r+i\m}$.\smallskip

We turn to the estimation of the Fourier coefficients of
$E_0^{\G(\qq)}(P^{\k_1},\n,i\m,g)$. We have seen that this suffices
for our purpose.

Each cusp $\k_1$ of $\G(\qq)$ has the form $-\frac{\delta}{\g}$, with
$\delta,\g \in \Okr$, such that the ideal $\fka=(\delta, \g)$ is
relatively prime to $\qq$, see for instance \cite{Gu}, Satz 1.

There are $\al, \bt \in \fka^{-1}$ such that $g_{\k_1}= \matc \dt
{-\bt} {-\g} \al$ sends $\infty$ to $\k_1= -\frac \dt\g$. If $\g_1$
runs through a set of representatives of
$\G(\qq)_{\k_1}\backslash\G(\qq)$, then $\g_1^{-1}\k_1$ runs through
the cuspidal orbit $\G(\qq) \k_1$. This can be reformulated in the
following way: the lower rows of $g_{\k_1}^{-1} \g_1 =\matc * * c d$
run through the pairs $(c,d) \in {\Okr}$ that satisfy
\begin{equation}\label{cuspconds}
\left \{
 \begin{split}
&\Okr c + \Okr d=\fka, \quad c\equiv \g \mod \qq \fka, \quad d\equiv
\dt \mod \qq \fka,\\
{}&\text{from each class }\vzm{(\e c,\e d)}{\e \in \Okr^\ast,\, \e
\equiv 1\mod \qq },\\ &\quad\text{ exactly one pair occurs}.
\end{split}
\right.
\end{equation}
 This follows from \cite{Gu}, Hilfsatz 1, Satz 1 (see also \cite{Ef},
Propositions 2.1 and 2.3).

As above, we put $g_{c,d}= g_{\k_1}^{-1}\g_1 =\matc * * c d$. Let
$z=x+iy \in \Hkr^d$. For $\re \n >\frac12$:
\begin{equation}
E_0^{\G(\qq)}(P^{\k_1}, \n,i\m,n[x]a[y]) = \sum_{(c,d)} a\left
[\frac{y}{|cz +d|^2}\right]^{\rho +2\n \rho + i\m}
\end{equation}
where $\left( \frac{y}{|cz +d|^2}\right)_j =\frac{y_j}{|c^{\s_j}z
+d^{\s_j}|^2}$ and where the pair $(c,d)$ runs over a set satisfying
the conditions in~\eqref{cuspconds}. By standard methods this gives
the following description of the Fourier term of order $r\ne 0$, at
the cusp $\infty$:
\begin{equation*}
\begin{split}
&\textstyle{\frac{1}{\vol ((\G(\qq)\cap N)\backslash N)}}
\int_{\G(\qq)_N\backslash N} \chi_r(n)^{-1} E_0^{\G(\qq)}(P^{\k_1},
\n,i\m,na[y])\,dn
\\ &= C_0 \sum_{(c,d)} |N(c)|^{-1-2\n}|c|^{-2i\m} e^{2 \pi i S(r d/c)}
|\pi r|^{\n + i\m} d_\infty^r (0,\n +i\m) W_{\infty,0}^{r,
\n+i\m}(a[y])
\end{split}
\end{equation*}
with $C_0$ a constant (see [BMP], (9),(16) for $d_\infty^r $ and
$W_{\infty,0}^{r, \n+i\m}$).

Note that
$$|c|^{-2\rho-4\n\r-2i\m} = \prod_j
\left|c^{\s_j}\right|^{-1-2\n-2i\m_j} = |N(c)|^{-1-2\n} |c|^{-2
i\m},$$ and that $c\mapsto |c|^{-2i\m}$ is a character on $F^\ast$,
vanishing on the group $\vzm{\e\in \Okr^\ast}{\e\equiv 1\bmod \qq}$,
since $\mu \in \Lkr_\k$.

Now the $(c,d) \in {\fka}\times{\fka}$ have to satisfy
\begin{equation}\label{conds}
\left\{
\begin{split}
& \Okr c + \Okr d = \fka,\quad c\equiv \g \text{ mod } \qq, \quad c\ne
0,\quad d\equiv \dt \text{ mod } \qq\\
& \text{there is exactly one representative of each class } d\text{
mod }\qq (c)\\
& \text{there is exactly one representative of each class}\\
&\quad c\vzm{\e \in \Okr^\ast}{\e\equiv 1 \text{ mod }\qq}.
\end{split}
\right.
\end{equation}

Since the ideal $\fka = \Okr \g + \Okr \dt$ is relatively prime to
$\qq$, the congruences $c\equiv \g ,\,d \equiv \dt \mod \qq$ are
equivalent to the same congruences$\mod \qq \fka$.

We shall denote by $\Psi(\g,\dt)$ a set of pairs $(c,d)$ satisfying
the conditions in (\ref{conds}) above.

In (10) and (18) of \cite{BMP}, we see that the quantity
$D_\k^{\infty,r}(\nu,i\m)$ in~\eqref{measure} is obtained by dividing
the Fourier term given above by the factors $d_\infty^r$ and
$W_{\infty,0}^{r, \n+i\m}$. We have thus arrived at a description of
 the following form
\begin{equation}\label{ddecomp}
D_\k^{\infty,r}(\nu,i\m)= \sum_{l=1}^{L} c_l a_l^{\n +i\m}
\Phi_r(\n+i\m; \g_l,\dt_l)
\end{equation}
with $c_l \in \CC,\, a_l \in (0,+\infty)^d$ and where the Dirichlet
series $\Phi_r$ is given for $\re \n> \frac12$ by:
\begin{equation}
\Phi_r(\n+i\m;\g,\dt) = \sum_{(c,d) \in \Psi(\g,\dt)}|N(c)|^{-1-2\n}
|c|^{-2i\m} e^{2\pi i S(rd/c)}.
\end{equation}

An estimate on the line $\re \n=0$ for the meromorphic continuation of
such a Dirichlet series will give an estimate of the same form for
$D_\k^{\infty,r}(\nu,i\m)$.

We will now consider one such series, depending on the pair $\g$ and
$\dt$ in $\Okr$ and where the ideal $\Okr \g + \Okr \dt=\fka$ is
coprime to $\qq$. One difficulty is the condition $\Okr c + \Okr
d=\fka$ in \eqref{conds} that prevents $d$ from running over a full
set of representatives. The method to approach this problem can be
found in \cite{Kloos}, \cite{Fr}, Ch. III \S 4, or in \cite{Ef}, \S
2.

Let $I_\qq$ be the group of ideals prime to $\qq$, and let $\m_\qq$ be
the M\"obius function on $I_\qq$. Let $\Psi_1\isdef\Psi_1 (\g,\dt)$
be the set of $(c,d) \in \fka\times \fka$ satisfying all the
conditions in (\ref{conds}), except for $\Okr c + \Okr d=\fka$. For
$(c,d) \in \Psi_1$, we only know that $\Okr c + \Okr d$ is an ideal
contained in $\fka$. Now
\begin{align*}
\Phi_r(\n+i\m;\g,\dt) =& \sum_{(c,d) \in \Psi_1}|N(c)|^{-1-2\n}
|c|^{-2i\m} e^{2\pi i S(rd/c)} \sum_{\fkb|\fka^{-1}(\Okr c + \Okr
d)}\m_\qq(\fkb)\\
=& \sum_{\fkb\subset \Okr : \fkb \in I_\qq} \m_\qq(\fkb) \sum_{(c,d)
\in \Psi_1: c,d \in \fka \fkb}|N(c)|^{-1-2\n} |c|^{-2i\m} e^{2\pi i
S(rd/c)}.
\end{align*}
Let $F_\qq =\vzm{(\al) \in I_\qq} { \al \in F^*, \al \equiv 1\mod \qq,
\al^{\s_j}>0, \text{ for } j=1,\dots, d}$. The quotient $I_\qq/F_\qq$
is the strict ray class group modulo $\qq$. We use this finite group
to split up $\Phi_r$ into the sum over $\tau \in I_\qq/F_\qq$ of
\begin{equation}\nonumber \Phi_r^\tau(\n+i\m;\g,\dt) =
\sum_{\fkb\subset \Okr : \fkb \in \tau} \m_\qq(\fkb) \sum_{(c,d) \in
\Psi_1: c,d \in \fka \fkb}|N(c)|^{-1-2\n} |c|^{-2i\m} e^{2\pi i
S(rd/c)}.
\end{equation}
Thus, to estimate $\Phi_r$ it suffices to estimate the finitely many
terms $\Phi_r^\tau$.

For a given $\tau$ we fix $\fkb_0 \in \tau$ an integral ideal. We
write each $\fkb \in \tau$ as $\fkb = \fkb_0(\vartheta)$, with
$\vartheta = \vartheta_\fkb \in F^*$ totally positive,
$\vartheta\equiv 1 \mod \qq$. For a given $\vartheta$ of this type,
we replace $c$ by $\vartheta c$ and $d$ by $\vartheta d$. The
conditions for the new pair $(c,d)$ are:
\begin{equation}\label{conds2}
\left\{\begin{array}{lll}
c\in \fka \fkb_0,\, c\ne 0, \quad c\equiv \g\bmod \qq\,\\ c \text{
modulo multiplication by } \e \in \Okr^\ast,\, \e\equiv 1 \mod \qq,\\
d \in \fka \fkb_0,\quad d\equiv \dt \mod \qq,\quad d \mod \qq
(c).\end{array}
\right.
\end{equation}
Let us denote by $\Psi_2(\fkb_0):=\Psi_2(\g,\dt;\fkb_0)$ a set of such
pairs $(c,d)$, as in \eqref{conds2}.

We have:
\begin{align}\nonumber \Phi_r^\tau(\n+i\m;\g,\dt)= \sum_{\fkb \subset
\Okr : \fkb \in \tau} \m_\qq(\fkb)&N(\fkb)^{-1-2\n} N(\fkb_0)^{1+2\n}
|\vartheta_\fkb|^{-2i\m}\\
\nonumber & \sum_{(c,d) \in \Psi_2(\fkb_0)}|N(c)|^{-1-2\n} |c|^{-2i\m}
e^{2\pi i S(rd/c)}.
\end{align}
We note that $\al \rightarrow |\al|^{2i\m}$ induces a unitary
 character of $F_\qq$. We fix a (unitary) extension $\ld_\m$ of this
character to $I_\qq$. For $\fkb=\fkb_0\, (\vartheta)$ as above we
have $|\vartheta|^{2i\m}=\ld_\m (\fkb)\ld_\m(\fkb_0)^{-1}$. So:
\begin{equation}\label{Phidecomp}
\Phi_r^\tau(\n+i\m;\g,\dt) = \ld_\m(\fkb_0)N(\fkb_0)^{1+2\n}\,
Q(\n,\ld_\m, \tau)\,\psi_r^{\fkb_0}(\n+i\m;\g,\dt),
\end{equation}
where
\begin{equation}\label{Q}
\begin{split}
Q(\n,\ld_\m;\tau)&=\sum_{\fkb\in \tau,\, \fkb\subset\Okr}\m_\qq(\fkb)
N(\fkb)^{-1-2\n} \ld_\m(\fkb)^{-1},\\
\psi_r^{\fkb_0}(\n+i\m;\g,\dt)&= \sum_{(c,d) \in
\Psi_2(\fkb_0)}|N(c)|^{-1-2\n} |c|^{-2i\m} e^{2\pi i S(rd/c)}.
\end{split}\end{equation}
We shall first estimate $\psi_r^{\fkb_0}(\n+i\m;\g,\dt)$. If we fix $c
\in \fka \fkb_0 , c\ne 0,\, c\equiv \g \mod \qq$, then the $d$ such
that $(c,d) \in \Psi_2(\fkb_0)$ run through representatives,$\mod
\qq(c)$, satisfying $d\equiv \dt'\mod \qq \fka \fkb_0$, for some
 $\dt'$ such that $\dt'\equiv \dt \mod \qq$ and $\dt'\equiv 0\mod \fka
\fkb_0$. Such $\dt'$ exists, since $\qq$ and $\fka \fkb_0$ are
relatively prime. So $d=\dt'+\vartheta$ with $\vartheta$ running
through all representatives of $\qq \fka \fkb_0 \mod \qq (c)$. Hence,
if $c$ is fixed and one sums over all $\vartheta$, one gets a non
zero contribution if and only if the character $x\rightarrow e^{2\p
iS(rx/c)}$ of $\qq \fka \fkb_0 \mod \qq (c)$ is trivial, that is, if
and only if $r\fkd\qq \fka \fkb_0\subset (c)$, where $\fkd$ is the
different of $\Okr$. So $\psi_r^{\fkb_0}(\n+i\m;\g,\dt)$ is given by
a sum over finitely many $c$, determined by
\begin{equation}\label{cconds}
\left\{\begin{array}{ll}\fka \fkb_0\,|\,(c)\\(c)\,|\,r\fkd\qq \fka
\fkb_0\\ c\equiv\g\mod \qq
\end{array}
\right.
\end{equation}
Hence
\begin{equation}\label{psiexpr}
\psi_r^{\fkb_0}(\n+i\m;\g,\dt)= \sum_{c \:{\rm as
\:in}\:\eqref{cconds} }|N(c)|^{-2\n} |c|^{-2i\m}N(\fka \fkb_0)^{-1}
\end{equation}
is a finite sum, hence everywhere holomorphic. For $\re \n =0$ it is
bounded by
\begin{equation}\label{psiest}
N(\fka \fkb_0)^{-1}\,\left| \vzm{\fkb\subset\Okr }{
\fkb\,|\,r\fkd\qq}\right|=\oh_\e\left(\left(|N(r)|N(\fkd)N(\qq)\right)^\e
\right),
\end{equation}
for each $\e>0$.\smallskip

We now turn to $Q(\n,\ld_\m;\tau)$, with $\tau \in I_\qq/F_\qq$. Set
\begin{equation}\label{Lfdef1}
L(s,\bar\ld_\m,\chi) = \sum_{\fkb \in I_\qq,\,
\fkb\subset\Okr}\frac{\overline{\ld_\m(\fkb)}
\chi(\fkb)}{N(\fkb)^{s}},
\end{equation}
where $\chi $ is a character of the ray-class group $I_\qq/F_\qq$.
Fourier analysis on $I_\qq/F_\qq$ allows to express
$Q(\nu,\ld_\m;\tau)$ in terms of these $L$-functions:
\begin{align}
\label{QL}
\left|I_\qq/F_\qq\right|\, Q(\nu,\ld_\m;\tau)
&= \sum_{\fkb\in I_\qq,\, \fkb\subset\Okr}\frac{
\overline{\ld_\m(\fkb)} \m_\qq(\fkb)}
{N(\fkb)^{1+2\nu}}\left(\sum_{\chi \in (I_\qq/F_\qq)^\wedge}
\overline{\chi(\tau)} \chi(\fkb)\right)\\
\nonumber \nonumber &=\sum_{\chi \in (I_\qq/F_\qq)^\wedge}
\frac{\overline{\chi (\tau)}}{L(1+2\nu, \bar\ld_\m,\chi)}.
\end{align}
where $(I_\qq/F_\qq)^\wedge$ denotes the character group of
$I_\qq/F_\qq$. So the task of estimating $Q(\nu,\ld_\m;\tau)$ for
$\re \nu=0$, is reduced to estimating $L(s, \bar\ld_\m,\chi)^{-1}$
for all $\chi \in (I_\qq/F_\qq)^\wedge$, for $\re s=1$.

In the case $F=\QQ$, there is a Riemann zeta function in the
denominator of the relevant expressions. There it suffices to use
$\frac1{\z(1+it)}=\oh(\log^7t)$ as $t\rightarrow\infty$, or the
better estimate $\frac1{\z(1+it)}=\oh(\log t)$, see \cite{Ti},
(3.6.5), and (3.11.8).

\subsection{Estimates for certain $L$-functions} In this subsection we
will obtain logarithmic estimates on the line $\re s=1$ for the
$L$-functions defined in~\eqref{Lfdef1}. In the simplest case,
$F=\QQ$, $\qq=\ZZ$, we can use the estimate
$\frac1{\z(1+it)}=\oh(\log^7t)$ as $t\rightarrow\infty$, or the
better estimate $\frac1{\z(1+it)}=\oh(\log t)$, see \cite{Ti},
(3.6.5), and (3.11.8). Such estimates are known for other
$L$-functions, like Dedekind zeta functions. We include a proof for
the general case we need, since we could not find it in the
literature in the needed generality.

We shall follow classical methods of Landau, leading to a
generalization of $\frac1{\z(1+it)}=\oh(\log^7t)$. Y.Motohashi kindly
pointed out to us a more involved method that should give a better
estimate, comparable to $\frac1{\z(1+it)}=\oh(\log t)$. In the
context of this paper we prefer to stick to an approach that is as
simple as possible.

As seen in the previous subsection, these estimates will allow us to
estimate the Fourier coefficients of Eisenstein series and therefore
to get a bound for the contribution of the Eisenstein term of the sum
formula.

We shall use the notations from~\S\ref{FcE}, but to keep this section
self-con\-tained, we will recapitulate the main ingredients.
\smallskip

Let $F$ be a totally real number field, $\qq$ an integral ideal in
$F$, $I_\qq$ the set of ideals in $\Okr=\Okr_F$ that are prime to
$\qq$ and $F_\qq=\vzm{(\al) \in I_\qq} {\al \equiv 1 \mod \qq,\, \al
\text{ totally positive}}$. So $I_\qq/F_\qq$ is the strict ray class
group modulo $\qq$. We consider $F^*$ embedded in $(\RR^*)^d$ in the
canonical way. Taking absolute values in each coordinate and denoting
by $V$ the composition of the absolute values with the embedding of
$F^*$, we obtain a map $V:\Okr^\ast\longrightarrow (\RR_{>0})^d$ with
a multiplicative lattice as its image in $(\RR_{>0})^d$ contained in
the subgroup defined by the equation $N(y)=1$. We have a unitary
character of $F^*/\vzm{\e\in \Okr^\ast}{\e\equiv1\bmod\qq}$ given by
$\ld_\m (x)= V(x)^{2i\m}$, where $y^{2i\m} =
\prod_{j=1}^dy_j^{2i\m_j}$. We have assumed that $ \sum_{j=1}^d \m_j
=0$ and we have extended the character $\ld_\m$ of $F_\qq$ to a
unitary character of the group~$I_\qq$. For simplification, we shall
denote this character by $\ld$ in the sequel. If $\m=0$, then we take
the extension $\ld$ equal to~$1$.

Our aim is to obtain a lower bound of the absolute value of
$L$-functions of the form
\begin{equation}\label{Lfdef}
L(s,\bar\ld,\chi) = \sum_{\fkb \in I_\qq,\,
\fkb\subset\Okr}\frac{\overline{\ld(\fkb)} \chi(\fkb)}{N(\fkb)^{s}},
\end{equation}
where $\chi$ is a unitary character of the group $I_\qq/F_\qq$, on the
line $\re s=1$. We follow the standard approach, which first
establishes upper bounds on lines $\re s=\s>1$. The latter bounds are
more easily obtained for ray class zeta functions:
\begin{equation}\label{Zdef} Z(s,\ld,\tau) \isdef \sum_{\fkb\in\tau,\,
\fkb\subset\Okr} \frac{\ld(\fkb)}{N(\fkb)^s} = \sum_{n=1}^\infty
\frac{f_\ld(n)}{n^s},\end{equation}
where $\tau$ is a class of $I_\qq/F_\qq$. The coefficients $f_\ld(n)$
of this Dirichlet series have less satisfactory estimates than the
partial sums
\begin{equation}\label{slddef} s_\ld(n) \isdef \sum_{k=1}^n f_\ld(k) =
\sum_{\fkb\in\tau,\, \fkb\subset\Okr,\, N\fkb \leq n} \ld(\fkb).
\end{equation}
In \S2, Chap.~IV of~\cite{Jz}, it is shown that
\begin{equation} \label{s1asympt}
s_1(n) = \al n + \oh(n^{1/d}),
\end{equation}
with $\al\neq0$ not depending on the ray class~$\tau$. We shall follow
the same method to prove Lemma~\ref{lem-s-est}, which shows that
$s_\ld(n) \ll n^{1/d}$ if $\ld\neq1$ for $d>1$. For $d=1$, we have
$s_1(n)=\al n + \oh(1)$.

It is more convenient to sum over numbers than over ideals. So we fix
 an integral ideal $\fkc$ in the inverse class $\tau^{-1}$, and write
$\fkb = (\al)\fkc^{-1}$ with $\al\in F^\ast$, $\al\equiv1\bmod\qq$,
$\al\gg0$ (i.e., $\al$ totally positive). We set
\begin{equation}
\label{sldpdef}
s'_\ld(n) \isdef \sum_\al\ld(\al) = \ld(\fkc) s_\ld(n/N(\fkc)),
\end{equation}
where $(\al)\fkc^{-1}$ runs through the integral ideals in $\tau$ such
that $N(\al) \leq n $. Since $|\ld({\mathfrak c})|=1$, the quantities
$s_\ld(n)$ and $s'_\ld(nN(\fkc))$ have the same growth behavior.

We proceed as in Lemma 2.5--2.7, Ch.~IV, in~\cite{Jz}. The $\al$ in
the sum satisfy
\begin{enumerate}
\item[i)] $\al\equiv\al_0 \bmod\qq\fkc$, where $\al_0\in\fkc$
satisfies $\al_0\equiv 1\bmod\qq$,
\item[ii)] $\al$ is totally positive,
\item[iii)] we use only one element from each class
$\vzm{\al\e}{\e\in\Okr^\ast,\, \e\gg0,\, \e\equiv1\bmod\qq}$,
\item[iv)] $N(\al)\leq n$.
\end{enumerate}
The $\al$ satisfying i) form a shifted lattice
$\Mkr=\al_0+\qq\fkc\subset F\subset\RR^d$. Condition~ii) requires
that we restrict the sum to $\al \in \Mkr \cap \RR^d_{>0}$.

We write each $x\in \RR^d_{>0}$ in the form $x=x_nx_u$, with $x_u$
such that $N(x_u)=1$, and $x_n=\left(
N(x)^{1/d},\cdots,N(x)^{1/d}\right)$. Let $\Fkr$ be a compact
fundamental domain for the multiplicative group $U_\qq\isdef
\vzm{\e\in\Okr^\ast}{\e\equiv1\bmod\qq,\, \e\gg1}$ in $N(x)=1$. Later
on we assume that the boundary of~$\Fkr$ is piecewise analytic. Let
$\Xkr_n$ be the set of all $x\in \RR^d_{>0}$ satisfying $N(x)\leq n$
and $x_u\in\Fkr$. The $\al$ in the sum may be chosen as the elements
of $\Mkr\cap \Xkr_n$.

For each totally positive $\al \in \Mkr$, we have $\ld(\al) =
e^{\sum_j 2i\m_j\log \al^{\s_j}}$. We extend this by defining $\ld(x)
= e^{\sum_j 2i\m_j\log x_j}$ for all $x\in \RR^d_{>0}$. We note that
$\ld=1$ on the multiplicative group~$U_\qq$, and $\ld(x)=\ld(x_u)$.
Thus

$$s'_\ld(n) = \sum_{x \in \Mkr \cap {\Xkr_n}} \ld (x).$$

Let $\Lambda$ be a fundamental domain for the lattice~$\qq\fkc$, of
which $\Mkr$ is a shift. We replace the sum $s'_\ld(n)$ by the
integral
\begin{equation} I_\ld(n) \isdef \frac1{\vol\Ld} \int_{\Xkr_n}\ld(y)\,
dy.
\end{equation}
The difference between $I_\ld(n)$ and $s'_\ld(n)$ consists of two
terms: The ``interior error'' occurs for each $x\in\Mkr$ such that
$x+\Ld \subset \Xkr_n$, whereas the boundary error is related to
those $x\in \Mkr$ for which $x+\Ld$ meets the boundary $\partial
\Xkr_n$ of~$\Xkr_n$.

We start with the interior error. For $x\in \Mkr$ such that
$x+\Ld\subset \Xkr_n$, we consider
\begin{align*}
\delta(x)\isdef& \ld(x) -\vol(\Ld)^{-1} \int_{x+\Ld}\ld(y)\,dy\\
=&\frac1{\vol(\Ld)} \int _\Ld (\ld(x) - \ld(x+\eta))\,d\eta
\end{align*}
For $\eta \in \Ld$, we have
$$\left(\frac{x_j+\eta_j}{x_j}\right)^{2i\m_j}= 1 + \oh_\Ld
(\m_j/x_j).$$
The set $\Fkr$ is compact, and has positive distance to all the
coordinate hyperplanes. So we find that $x_j \ge C_1 N(x)^{1/d}$, for
some $C_1 >0$. We have:
\begin{align*}
\left(\frac{x+\eta}{x}\right)^{i\m}&=\prod_{j=1}^d \left(1 +
\oh_{\Ld,\Fkr}(|\m_j| N(x)^{-1/d} )\right)\\
&= 1 +\oh_{\Ld,\Fkr}(\parallel
\m
\parallel N(x)^{-1/d})
\end{align*}
with $\parallel \m \parallel=\max_{j=1}^d |\m_j|$. Thus,
\begin{align*}
\delta(x)&= \frac{\ld(x)}{\vol (\Ld)} \int _\Ld
\oh_{\Ld,\Fkr}\left(\parallel \m
\parallel N(x)^{-1/d} \right)\, d\eta \ll_{\Ld, \Fkr}
\parallel\m\parallel N(x)^{-1/d},
\end{align*}
since $|\ld (x)|=1$. We find the following estimate for the interior
error:
\begin{equation}
\sum_{x \in \Mkr \cap \Xkr_n}|\dt (x)|\ll_{\Ld,\Fkr}
\parallel\m\parallel \sum_{x \in \Mkr \cap \Xkr_n} N(x)^{-1/d}
\end{equation}

We use partial summation to handle the latter sum on the basis of the
estimate $s_1'(n) \ll n$ for the case $\ld=1$, (see Statement~2.15,
Chap.~IV of~\cite{Jz}):
\begin{align*}
&\sum_{x\in \Mkr \cap \Xkr_n} N(x)^{-1/d} \le \sum_{m=1}^{n} m^{-1/d}(
s_1'(m) - s_1'(m-1)) \\
&\qquad \le \sum_{m=1}^{n-1} s_1'(m)\left( m^{-1/d}-
(m+1)^{-1/d}\right) + n^{-1/d} s_1'(n)\\
&\qquad \ll \sum_{m=1}^{n-1} m^{-1/d}+ n^{-1/d+1} \ll n^{1-1/d}
\end{align*}
Thus
\begin{equation}\label{interr}
\sum_{x\in \Mkr \cap \Xkr_n^a}\dt(x)=\oh_{\Ld,\Fkr}(\parallel \m
\parallel n^{1-1/d})
\end{equation}
\smallskip

At the boundary $\partial\Xkr_n$, for some $x\in \Mkr \cap \Xkr_n$,
the set $x+\Ld$ may jut out of $\Xkr_n$. At other points we may have
for $x \in \Mkr \setminus \Xkr_n$ that $(x+\Ld)\cap \Xkr_n\ne
\emptyset$. We note that $\Xkr_n$ is obtained from $\Xkr_1$ via
multiplication by $t=n^{1/d}$. We apply the argument in the proof of
Theorem~2, \S2, Chap.~VI, in \cite{La68}, to see that the number of
$x\in \Mkr$ such that $(x+\Ld)\cap \partial\Xkr_n\ne \emptyset$ is
\begin{equation}\label{bderr}
\oh(n^{1-1/d})
\end{equation}
The boundary error has the same estimate (not depending on~$\m$).
\smallskip

We now turn to the integral $I_\ld(n)$. It is convenient to change
coordinates to $z_j=\log\,y_j.$ Then the region of integration
corresponds to the infinite region
\begin{equation}
\Ykr_n=\left\{\begin{array}{ll} S(z) \le \log n\qquad
S(z)=\sum_{j=1}^d z_j,\\
z_u \in \Hkr \qquad z_u= z-(\frac{S(z)}{n},\dots,\frac{S(z)}{n}),
\end{array}
\right.
\end{equation}
where $\Hkr\isdef \log\Fkr$. We note that $\Hkr$ is a fundamental
domain for a group of translations in $S=0$ for which $z\mapsto
e^{2i\m.z}$ is periodic. Therefore:
\begin{equation} \label{Inull}
I_n = \frac{1}{\vol\Ld} \int_{\Ykr_n}e^{2i\m \cdot z}e^{S(z)}\,dz.
\end{equation}
The factor $e^{S(z)}$ and the restriction $z_u\in \Hkr$ are
responsible for the absolute convergence. If we integrate first over
a hyperplane $S(z)=C$, $C$ a constant, we get zero; therefore
$I_n=0$. In the light of \eqref{sldpdef}, \eqref{interr},
\eqref{bderr} and~\eqref{Inull}, we get:
\begin{lem}\label{lem-s-est}If $\ld$ is non-trivial, then the sum
$s_\ld(n)$ defined in~\eqref{slddef} satisfies
\begin{equation}\label{sldest}
s_\ld (n)=\oh_{\Ld,\Fkr,\fkc}(\parallel \m \parallel n^{1-1/d} ).
\end{equation}
\end{lem}
\medskip

Now we estimate the zeta function $Z(s,\ld,\tau)$, see~\eqref{Zdef}.
The argument will be essentially the same as in~\cite{La}, and will
only be sketched. We will follow the standard Mertens' scheme, using
\eqref{sldest} to keep track of the influence of $\ld=\ld_\m$. For
$\re s>1$, one finds with \eqref{slddef}:
\begin{align*}
Z(s,\ld,\tau)= \sum_{n=1}^m \frac{f_\ld(n)}{n^s}+ \sum_{n=m+1}^\infty
s_\ld(n)(n^{-s} - (n+1)^{-s}) -\frac{s_\ld(m)}{(m+1)^s},
\end{align*}
which stays valid for $\re s>1-\frac1d$ if $\ld\neq1$. One takes $\ld
\ne 1$, $\s =\re s\geq1$ and uses \eqref{sldest} to obtain:
\begin{align*}
\, &\left|\sum_{n=1}^m \frac{f_\ld(n)}{n^s}\right| \le
\sum_{n=1}^{m-1} s_1(n) \left(\frac1n -\frac1{n+1}\right) +
\frac{s_1(m)}m \ll \log m \smallskip\\
\, & \left|\sum_{n=m+1}^\infty s_\ld(n)(n^{-s} - (n+1)^{-s})\right|
\ll |s|\parallel \m \parallel m^{-1/d} \smallskip\\
\,& \left|\frac{s_\ld(m)}{(m+1)^s}\right| \ll \parallel \m
\parallel m^{-1/d}
\end{align*}
A suitable choice is $m=[1+(1+|s|)^d\parallel \m
\parallel^d]$.
This gives
\begin{equation}
Z(s,\ld,\tau)\ll_d \log(1+ |s|) + \log \parallel \m \parallel,
\end{equation}
 for $\ld \ne 1$, $\re s\geq 1$.

Next, one estimates the derivative of the zeta function:
\[ Z'(s,\ld,\tau) = - \sum_{n=1}^\infty \frac{f_\ld(n)\log n}{n^s}.\]
An argument entirely analogous to the one above gives for $\re s\geq
1$, $\ld\neq1$:
\begin{equation}
Z'(s,\ld,\tau)\ll \log^2 (1+|s|) + \log^2 \parallel \m \parallel.
\end{equation}

The classical approach, c.f.\ Landau, \cite{La} p.~91, p.~94, implies
similar estimates for the case $\ld=1$, $1\leq \re s\leq 2$, with
bounds $\frac1{|\im s|^l} + \log^l|\im s|$, $l=1$ or~$2$, for $Z$ and
$Z'$ respectively, instead of $\log^l(1+|s|) + \log^l \parallel \m
\parallel$.
\medskip

We now turn to the $L$-series defined in~\eqref{Lfdef}. It is a linear
combination of series $Z(s,\ld,\tau)$, with $\tau$ running through
the ray classes. The estimates for $Z$ imply similar estimates for
$L$. {}From~\eqref{s1asympt} it follows that $Z(s,1,\tau) =
\frac\al{s-1} + h_\tau(s)$, with $h_\tau$ holomorphic on $\re
s>1-\frac1d$. If the character $\ch$ of the ray class group
$I_\qq/F_\qq$ is non-trivial, we have $L(s,1,\ch) = \sum_\tau
\ch(\tau) h_\tau(s)$, holomorphic on $\re s>1-\frac1d$. So the
finitely many $L(s,1,\ch)$ with $\ch\neq1$ have a common upper bound
for $1\leq \re s \leq 2$, $|\im s|\leq 1$, and so do their
derivatives. Of course, $L(s,1,1)$ has a first order singularity at
$s=1$. We summarize:

\begin{alignat*}2 L(s,\ld,\ch) &\ll \log(2+|t|) + \z_\ld\log
\parallel\m\parallel + \dt_{\ld,\ch} |t|^{-1} &=:& \, b(t,\m,\ch),\\
L'(s,\ld,\chi) &\ll \log^2(2+|t|) + \z_\ld\log^2
\parallel\m\parallel + \dt_{\ld,\ch} |t|^{-2} &\ll& b(t,\m,\ch)^2,
\end{alignat*}
for $s=\s+it$, $\s\geq1$, $t\neq0$, $\ch$ a character of the ray class
group $I_\qq/F_\qq$. We use $\dt_{\ld,\ch}=1$, if $\ch=1$, $\ld=1$,
and $\dt_{\ch,\ld}=0$ otherwise, and $\z_\ld=1$ if $\ld\neq 1$,
$\z_1=0$. In the case $\dt_{\ld,\ch}=0$, these estimates extend to
$s=1$ by continuity.

The product expansion of the $L$-function implies
\[ \log |L(\s+it,\ld,\chi)|= \sum _{n=1}^\infty \frac{1}{n}
\sum_{{\fkp} : (\fkp,\qq)=1}{N(\fkp)}^{-n\s} \re \left(\ld(\fkp)^n
\chi(\fkp)^n {N(\fkp)}^{-nit}\right) \]
{}From the inequality $3+4\cos\phi +\cos 2\phi \ge 0$ it follows in
the standard way that for $\s>1$
\[ |L(\s,1,1)|^3|L(\s+it,\ld,\chi)|^4 |L(\s +2it,\ld^2,\chi^2)| \ge 1.
\]
We suppose that $t\neq0$ or $\dt_{\ld,\ch}\neq0$. Then we have for
$\s>1$:
\[|L(\s+it,\ld,\chi)|\ge \frac{ C_2'
(\s-1)^{3/4}}{|L(\s+2it,\ld^2,\ch^2)|^{1/4}}.\]
leaving the exceptional case $\m=0$, $\ch^2=1$, $|t|\leq1$ aside for
the moment, we have $b(2t,2\m,\ch^2) \ll b(t,\m,\ch)$, and
obtain\begin{equation}
|L(\s+it,\ld,\chi)|\ge C_2 \frac{(\s-1)^{3/4}}{b(t,\m,\ch)^{1/4}}.
\end{equation}
Furthermore,
\begin{align*}
&\left|L(\s+it,\ld,\chi)-L(1+it,\ld,\chi)\right|
\\
&\quad\ll (\s-1)|L'(1+\theta +it,\ld,\chi)|,
\quad 0< \theta <\s-1\\
&\quad\ll (\s-1) b(t,\m,\ch)^2.
\end{align*}
For $t\neq0$ if $\dt_{\ld,\ch}=0$:
\begin{align*}
&|L(1+it,\ld,\chi)|\ge
|L(\s+it,\ld,\chi)|-\left|L(\s+it,\ld,\chi)-L(1+it,\ld,\chi)\right|\\
&\quad \ge \frac{C_2 (\s-1)^{3/4}}{b(t,\m,\ch)^{1/4}} -C_3 (\s-1)
b(t,\m ,\ch)^{2}.
\end{align*}
One may choose
\[ \s=1+\left(\frac{C_2}{2C_3}\right)^4 b(t,\m,\ch)^{-9},\qquad
C_4=\left(\frac{C_2}{2C_3}\right)^3\frac{C_2}{2}, \]
to obtain $|L(1+it,\ld,\ch)|\geq C_4\, b(t,\m,\ch)^{-7}$.

In the exceptional cases when $\m=0$ (for which we have chosen
$\ld=1$), $\ch^2=1$, $|t|\leq 1$, we find
\[ |L(\s+it,1,\ch)| \gg (\s-1)^{3/4} |\s+2it-1|^{1/4}.\]

So these $L(s,1,\ch)$ cannot have a zero at $\s\geq 1$, $0<|t|\leq 1$. 
We know that there is a pole at $s=1$ if $\ch=1$ (see, e.g.,
\cite{Jz}, Ch.~IV, Prop.~4.1). If $\ch\neq1$, then $L(s,1,\ch)$ is
holomorphic and non-zero at $s=1$, see~\cite{Jz}, Ch.~V, Prop.~10.2.
Thus, $L(s,1,\ch)^{-1}$ is bounded on the region $\s\geq1$, $|t|\leq
1$. The ray class group $I_\qq/F_\qq$ is finite, so there are only
finitely many exceptional~$\ch$.

We have thus proved:
\begin{prop}\label{lemLest}Let $\ld$ be a unitary character of $I_\qq$, which is
 given by $\ld(\x) =\prod_{j=1}^d |\x^{\s_j}|^{2i\m_j}$ on the
ideals in $I_\qq$ of the form $(\x)$, $\x\in F^\ast$, and which is
equal to~$1$ if~$\m=0$. Let $\ch$ be a character of the ray class
group $I_\qq/F_\qq$. Then, for $\re s =1$, we have:
\[ \frac1{|L(s,\ld,\chi)|}\ll
\begin{cases} \log^7(2+|\im s|) + \log^7 \|\m\| &\text{ if }
\ld\neq1,\\
\log^7 (2+|\im s|) &\text{ if } \ld=1.
\end{cases}
\]
\end{prop}

In view of \eqref{QL}, the function $Q(\n,\ld,\tau)$ satisfies the
same estimates:
\begin{equation} \label{Qest}
Q(\n,\ld,\tau) \ll \log^7(2+|\im \n|) + \z_\ld \log^7 \|\m\|
\end{equation}
for $\re \n=0$, with $\z_\ld=1$ if $\ld\neq1$, $\z_1=0$.

\subsection{Estimation of the term corresponding to the continuous
spectrum} We turn to the estimation of the term $\int_Y k(\n)\,
d\s_{r,r}^{\rm cont}(\n)$ in~\eqref{measure}, for the test function
$k = \times_j k_j$ with arbitrary $k_j\in \Lkr$ if $j\in E$, and
$k_j$ as indicated in \eqref{Gpchoice}, \eqref{Gmchoice} for the
other places. {}From \eqref{ddecomp}, \eqref{Phidecomp},
\eqref{psiexpr}, \eqref{psiest}, and \eqref{Qest}, we conclude for
each cusp~$\k$:
\begin{align*}
D^{\infty,r}_\k(iy,i\m) &\ll |N(r)|^\e \left(\log^7(2+|y|) + \log^7
\|\m\| \right)\qquad\text{ if }\m\neq0,\displaybreak[0]\\
 D^{\infty,r}_\k(iy,0) &\ll |N(r)|^\e \log^7(2+|y|).
\end{align*}
In view of~\eqref{measure} and~\eqref{Gmchoice}, the integral vanishes
if $Q_-\neq0$. So we need only look at the case $\vz{1,\ldots,d} = E
\sqcup Q_+$.
\begin{align*}
\int_Y k \, &d\s_{r,r}^{\rm cont} \ll |N(r)|^\e \sum_\k
c_\k\sum_{\m\in \Lkr_\k} \int_{-\infty}^\infty \|k\|_{\al, a,E}
\prod_{j\in E} \left( 1+|y+\m_j|\right)^{-a} \displaybreak[0]\\
&\qquad\hbox{} \cdot \prod_{j\in Q_+} e^{-s(y+\m_j)^2-s/4} \cdot
\left( \log^7(2+|y|) + \z_\m \log^7 \|\m\| \right)\, dy,
\end{align*}
with $\z_0=0$, $\z_\m=1$ if $\m\neq0$. For each of the finitely many
$\k\in \Pkr$, we estimate the sum over $\Lkr_\k$ by an integral over
the hyperplane $\sum_j x_j=0$. For $\al\in
\left(\frac12,\tau\right]$:
\begin{align}\nonumber \int_Y k &\, d\s_{r,r}^{\rm
cont}\displaybreak[0]\\
\nonumber
&\ll |N(r)|^\e \|k\|_{\al,a,E} \int_{\RR^d} \prod_{j\in
E}\left(1+|x_j|\right)^{-a} \prod_{j\in Q_+} e^{-s x_j^2}
\displaybreak[0]\\
\nonumber
&\qquad\hbox{} \cdot {\textstyle \left( \log^7\left(
2+\left|{\textstyle\sum_j x_j} \right| \right) + \log^7\left( 1 +
\max_j \left|{\textstyle x_j - d^{-1}\sum_\ell x_\ell} \right|
\right) \right)}\, dx\displaybreak[0]\\
\nonumber
&\ll |N(r)|^\e \|k\|_{\al,a,E} \prod_{j\in E} \int_{-\infty}^\infty
\left( 1+ |x| \right)^{-a+\e} \, dx\ \displaybreak[0] \\
\nonumber
&\qquad\hbox{} \cdot \prod_{j\in Q_+} \int_{-\infty}^\infty (1+|x|)^\e
e^{-sx^2}\, dx \displaybreak[0] \\
\label{ETest}
&\ll |N(r)|^\e \|k\|_{\al,a,E} s^{(|E|-d)(\frac {1+\e}2)}, \text{ if }
Q_-=\emptyset.
\end{align}
On the other hand, we have seen that $\int_Y k \, d\s_{r,r}^{\rm
cont}= 0$ , if $Q_-\ne\emptyset$.

\section{Proof of Theorem~\ref{propv}}\label{discrterm} The last term
to be considered is given by the measure $d\s_{r,r}^{\rm discr} =
d\s_{r,r}-d\s_{r,r}^{\rm cont}$:
\begin{equation}\label{dmeasure}
 \int_Y f(\n)\, d\s_{r,r}^{\rm discr} (\n) = \sum_{\irr\neq \one}
f(\n_\irr) \left| c^r(\irr)\right|^2.\end{equation}
The sum formula~\eqref{sf*}, and the estimates in Proposition
\ref{Dtest1}, Proposition \ref{KLTest}, and \eqref{ETest}, give for
test functions $k$ as chosen in \S\ref{prfpropv}, and
$\frac12<\al\leq\tau$, $\e>0$ sufficiently small:
\begin{align}
\label{discrest}
\int_Y &k(\n) \, d\s_{r,r}^{\rm discr}(\n) =
\tfrac{2^{1+|E|}}{(2\p)^d} \sqrt{\left| D_F\right|} \Eta_E(k)
s^{|E|-d} \displaybreak[0]\\
\nonumber
&\hbox{} +\oh_{F} \left( \|k\|_{\al,a,E} \; s^{|E|-d+c}\right)
\displaybreak[0]\\
\nonumber
&\hbox{}+ \text{(if $Q_-=\emptyset$)\ } \oh_{F, \qq,\e} \left(
\|k\|_{\al,a,E}|N(r)|^\e s^{(|E|-d)(\frac{1+\e}2)} \right)
\displaybreak[0] \\
\nonumber
&\hbox{} + \oh_{F,\qq,\e}\left( \|k\|_{\al,a,E}
N(r)^{1/2+\e}s^{-(3/4+\e)|Q_+|-(1/4+1/8\al+\e) |Q_-|} \right).
\end{align}

Taking into account that
\begin{align} \int_y k(\n)\, d\s_{r,r}^{\rm discr} (\n) =
\sum_{\substack{\irr\neq\one\\
\ld_{\irr,j}\geq0,\, j\in Q_+\\
\ld_{\irr,j}<0,\, j\in Q_- }} \left| c^r(\irr)\right|^2 e^{-s{\|
\ld_{\irr,Q} \|}_1 } \prod_{j\in E} k_j(\n_{\irr,j}),
\end{align}
and that $d=|E|+|Q_+|+|Q_-|$, we obtain the following result:
\begin{prop}\label{errorZeta} Let $r\in\Okr'\setminus\vz0$. Choose a
partition $E$, $Q_+$, $Q_-$ of $\vz{1,\ldots,d}$ with $Q =Q_+ \cup
Q_- \ne \emptyset$ and let $\Rkr$ be the corresponding set as
in~\eqref{partition}. Let $g=\times_{j\in E}\, g_j$, with $g_j\in
\tilde\Lkr$ for each $j\in E$, let $\frac12<\al\leq \tau<1$ and
$\e>0$ sufficiently small. If $Z_s(g)$ is as in \eqref{Zeta}, then we
have
\begin{equation*}
\begin{split}
 &Z_s(g) ={\textstyle\frac{ 2^{1+|E|}}{(2\pi)^d}}\sqrt{|D_F|}\,
\prod_{j\in E} \tilde\Eta(g_j)s^{-|Q|} + \oh_F (\|k\|_{\al,a,E}
s^{-|Q| +c}) \displaybreak[0] \\
&+\;\hbox{(if $Q_-=\emptyset$)\ }\oh\left(\|k\|_{\al,a,E} |N(r)|^\e
s^{-(1/2+\e)|Q|} \right)\\
&+\oh_{F,\qq,\e}\left( \|k\|_{\al,a,E}
N(r)^{1/2+\e}s^{-\left((3/4+\e)|Q_+|\,+\,(1/4+1/8\al+\e)|Q_-|\right)}
\right).
\end{split}
\end{equation*}
Here $c=\frac 12$ if $Q_-\ne \emptyset$ and $c=1$ otherwise, and
$\tilde\Eta(g)$ is as given in \eqref{Etatil}. Also, if $E=\emptyset$
then the product factor in the right-hand side should be interpreted
as $1$.
\end{prop}
We compare the three error terms in the $s$-aspect, under the
assumptions $\e<\frac1{20}$ and $\e<\frac14-\frac1{8\al}$.

If $|Q_-|\geq 1$, then
$-\left(\frac34+\e\right)|Q_+|-\left(\frac14+\frac1{8\al}+\e\right)|Q_-|
> -|Q_+|-|Q_-|+\frac12|Q_-| \geq -|Q_+| -|Q_-|+\frac12$.

If $|Q_-|=0$, then $-\left(\frac12+\e\right)|Q_+| \geq
-\left(\frac34+\e\right)|Q_+|$. It turns out that $-|Q_+|+1$ is
smaller than $-\left(\frac34+\e\right)|Q_+|$ precisely if $|Q_+|\geq
5$. Hence the error term in Proposition~\ref{errorZeta} is
\begin{alignat*}2
&\oh_{F,k,E,r,\qq}\left( s^{-|Q|+1/2} \right)&&\text{ if
}|Q_-|\geq1,\\
&\oh_{F,k,E,r,\qq}\left( s^{-(3/4+\e)|Q|} \right)&&\text{ if
}|Q_+|\leq 4,\, |Q_-|=0,\\
&\oh_{F,k,E,r,\qq}\left( s^{-|Q|+1} \right)&&\text{ if }|Q_+|\geq 5,\,
|Q_-|=0.
\end{alignat*}\smallskip

Taking into account the relation \eqref{g} between the test functions
$k$ and $g$, we obtain Theorem~\ref{propv}.

\newcommand\bibit[5]{
\bibitem[#1]
{#2}#3: {\em #4;\/ } #5}

\end{document}